\documentclass[twocolumn]{autart}
\usepackage[T1]{fontenc}
\usepackage[latin9]{inputenc}
\usepackage{color}
\usepackage{enumerate,graphicx,epstopdf}
\usepackage{amsmath,amssymb,bm}
\usepackage{cite}
\usepackage{mathtools}
\usepackage{array}
\usepackage{algorithm,algorithmic}
\usepackage{color}
\usepackage{lipsum}

\newtheorem{ass}{Assumption}
\newtheorem{rmk}{Remark}

\newtheorem{dfn}{Definition}
\newtheorem{lmm}{Lemma}

\newcommand{\reals}{\mathbb{R}}
\newcommand{\nnints}{\mathbb{Z}_+}

\newcommand\veps[0]{\varepsilon}
\newcommand\eps[0]{\epsilon}
\newcommand\Kinf[0]{\mathcal{K}_\infty}
\newcommand\KL[0]{\mathcal{KL}}
\newcommand\K[0]{\mathcal{K}}
\newcommand\limsupk[0]{\underset{k\to\infty}{\overline{\lim}}}
\newcommand\id[0]{\mathrm{id}}
\newcommand\mc[1]{\mathcal{#1}}
\renewcommand\int[0]{\mathrm{Int~}}
\newcommand\Oh[0]{\mathcal{O}}
\newcommand\rad[1]{\mathrm{rad~#1}}
\newcommand\ints[0]{\mathbb{Z}}
\newcommand\dom[0]{\mathrm{dom}~}

\begin{document}

\begin{frontmatter}

\title{Time Distributed Optimization for Model Predictive Control: Stability, Robustness, and Constraint Satisfaction\thanksref{footnoteinfo}}

\thanks[footnoteinfo]{This research is supported by the Toyota Research Instituite (TRI) and by the National Science Foundation through awards  CMMI 1904441 and CMMI 1562209. TRI provided funds to assist the authors with their research but this article solely reflects the opinions and conclusions of its authors and not TRI or any other Toyota entity.}

\author[UM]{Dominic Liao-McPherson}\ead{dliaomcp@umich.edu},
\author[UCB]{Marco M. Nicotra}\ead{marco.nicotra@colorado.edu},
\author[UM]{Ilya Kolmanovsky}\ead{ilya@umich.edu}

\address[UM]{Department of Aerospace Engineering, University of Michigan, 1221 Beal Avenue, Ann Arbor, MI 48109 } 
\address[UCB]{Department of Electrical, Computer, and Energy Engineering, University of Colorado Boulder, 425 UCB, Boulder, CO 80309}


\begin{abstract}
Time distributed optimization is an implementation strategy that can significantly reduce the computational burden of model predictive control by exploiting its robustness to incomplete optimization. When using this strategy, optimization iterations are distributed over time by maintaining a running solution estimate for the optimal control problem and updating it at each sampling instant. The resulting controller can be viewed as a dynamic compensator which is placed in closed-loop with the plant. This paper presents a general systems theoretic analysis framework for time distributed optimization. The coupled plant-optimizer system is analyzed using input-to-state stability concepts and sufficient conditions for stability and constraint satisfaction are derived. When applied to time distributed sequential quadratic programming, the framework significantly extends the existing theoretical analysis for the real-time iteration scheme. Numerical simulations are presented that demonstrate the effectiveness of the scheme.
\end{abstract}

\end{frontmatter}

\section{Introduction} \label{ss:intro}
Model Predictive Control (MPC) \cite{grune2017nonlinear} is a widely used control technique that computes control actions by solving an Optimal Control Problem (OCP) over a finite receding horizon. MPC can systematically handle constraints and nonlinearities but is challenging to implement since it requires the solution of a constrained and potentially non-convex OCP at each sampling instant. The development of robust and efficient quadratic and convex programming solvers, see e.g., \cite{liaomcphersonFBRS,patrinos2014accelerated,domahidi2013ecos,rao1998application}, has enabled the application of linear-quadratic MPC to a wide variety of systems. However, the implementation of MPC for systems with limited onboard computing power, fast sampling rates, and/or pronounced nonlinear dynamics remains an open problem.

\begin{figure}[h]
	\centering
	\includegraphics[width=0.95\columnwidth]{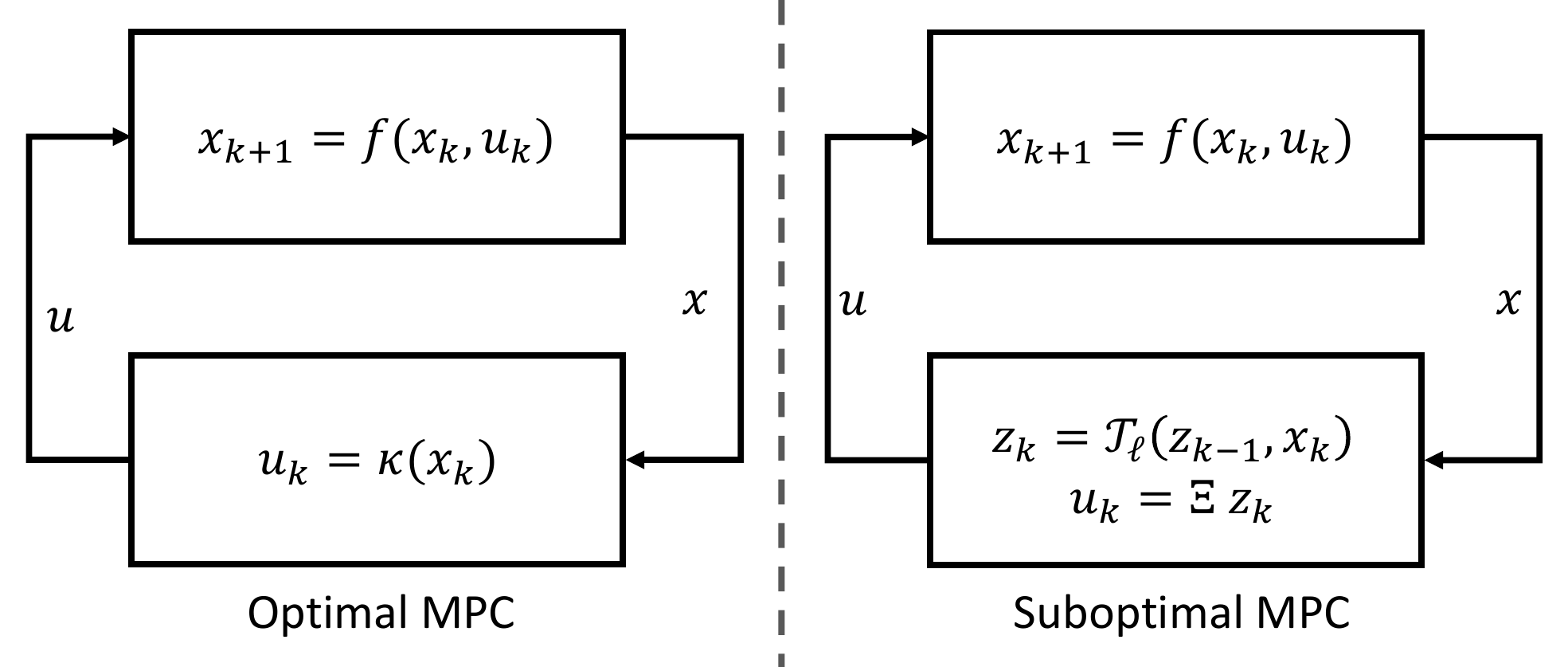}
	\caption{A comparison of suboptimal MPC with TDO and optimal MPC. The $\kappa$ operator represents the optimal MPC feedback law, the $\mathcal{T}_\ell$ operator represents $\ell$ iterations of an optimization algorithm, and $\Xi$ is a selection matrix that extracts the control action. One can roughly identify $\kappa(x) = \Xi \mc{T}_\infty(z_0,x)$ for any $z_0$ for which $\mc{T}$ converges.}
	\label{fig:subopt_mpc_fig}
\end{figure}

One approach for reducing the computational cost of MPC is time distributed optimization (TDO). TDO distributes optimizer iterations over time by exploiting the robustness of MPC to suboptimality \cite{ALLAN201768,pannocchia2011conditions,limon2009input}. Rather than accurately solving the OCP at each sampling instant, TDO maintains a guess of the optimal solution and improves it at each timestep by performing a finite number of iterations of an optimization algorithm. TDO can be interpreted as a dynamic compensator that maintains a solution estimate as an internal state, the dynamics of which are defined by the optimizer iterations. As illustrated in Figure~\ref{fig:subopt_mpc_fig}, this interpretation differs from ``ideal'', or ``optimal'' MPC which is an implicitly defined static feedback law.

There are a variety of TDO variants proposed in the literature. The stability of input constrained TDO controllers using linearly convergent optimization algorithms are studied in \cite{graichen2010stability}. Unconstrained suboptimal NMPC without terminal conditions is considered in \cite{grune2010analysis}. A fixed point scheme for input constrained MPC of sampled data input affine systems is proposed in \cite{graichen2012fixed}, a gradient based dynamic programming approach is considered in \cite{steinboeck2017design}, a proximal gradient method for linear input-constrained MPC is studied in \cite{van2019real}, and a continuous time gradient flow based approach is described in \cite{liao2018embedding}. These methods use some combination of shifting terminal control updates and first order optimization methods. In \cite{ALLAN201768,pannocchia2011conditions}, a generic suboptimal MPC scheme is considered and sufficient conditions on the warmstart for robust stability are derived; the optimization algorithm is not specified and its convergence is not considered. The robustness of MPC to disturbances arising from incomplete optimizations is considered in \cite{zavala2009advanced,gros2016linear} and conditions for complexity certification of suboptimal state constrained linear MPC are presented in \cite{rubagotti2014stabilizing}. However, the treatment of the optimizer itself as a dynamic system was not pursued.

An alternative to gradient based approaches are second order methods. In particular, Time Distributed Sequential Quadratic Programming (TD-SQP) methods are attractive since they can be implemented using existing Quadratic Programming (QP) solvers. The fundamental idea behind a TD-SQP based model predictive controller is to apply a finite number of SQP iterations at each sampling instant and to \textit{warmstart} the iterations with the solution estimate from the previous sampling instant. A widely used variant of TD-SQP is the Real-Time Iteration (RTI) scheme \cite{diehl2005real} which uses a Gauss-Newton Hessian approximation and performs a single SQP iteration per sampling instant. The RTI scheme has been successfully applied to a variety of applications including engines \cite{albin2018vehicle,zhu2017economic}, kites\cite{ilzhofer2007nonlinear}, cranes \cite{vukov2012experimental}, ground vehicles \cite{frasch2013auto}, race cars\cite{liniger2015optimization}, distillation columns \cite{diehl2001real} and wind turbines \cite{gros2013economic}. Software for implementing the RTI scheme is provided by the ACADO toolkit \cite{houska2011acado}. Despite its widespread success, formal stability guarantees for the RTI scheme have only been provided in the absence of inequality constraints \cite{diehl2005nominal}.

It should be noted that TDO is distinct from so-called suboptimal solution tracking, sensitivity, or running methods, e.g., \cite{zavala2009advanced,zavala2010real,liaomcpherson2018,dinh2012adjoint,ghaemi2009integrated,jaschke2014fast} which are tailored numerical methods for tracking the solutions of parameterized nonlinear programs/generalized equations. These methods are typically used to accelerate or replace existing nonlinear programming solvers to reduce computation times, which is different from considering the dynamic interactions between the plant and the optimizer. Some of them, e.g., \cite{zavala2009advanced}, consider robust stability by treating suboptimality as a bounded disturbance. This differs from our approach where we treat suboptimality as the output of a dynamic system which is coupled with the closed-loop plant.

This paper begins by presenting a system theoretic framework for analyzing a broad class of TDO algorithms. Specifically, the framework applies to any MPC feedback law that is Locally Input-to-State Stable (LISS) combined with any optimization algorithm featuring a convergence rate that is at least locally $q$-linear. Any MPC formulation with proven LISS properties can be used, including nominal MPC with terminal constraints \cite{limon2009input}, robust MPC\footnote{A \textit{nominal} MPC controller does not explicitly consider the presence of disturbances in the OCP formulation, unlike a \textit{robust} MPC controller.} \cite{limon2009input}, and MPC formulations with no terminal constraints \cite[Chapter 6]{grune2017nonlinear}. in this paper, we establish the existence of a joint region of attraction for the state and solution estimate, i.e., we show that if the initial state is sufficiently close to the origin and if the initial solution estimate is sufficiently accurate, the state will converge to the origin and the estimate will converge to the optimal solution. Moreover, we analyze the effect of performing more iterations, establish robustness properties, and show that, if the initial solution guess is within the convergence basin of the optimization method, TDO can recover the robust region of attraction of optimal MPC with a finite number of iterations.

The proposed theoretical framework is then specialized to the RTI scheme. 
Our analysis extends that in \cite{diehl2005nominal} as follows:
(i) We explicitly consider inequality constraints and relax the terminal state constraint to a terminal set constraint;
(ii) We explicitly consider the robustness properties of the RTI scheme by establishing LISS of the closed-loop system;
(iii) We analyze the effect of the number of SQP iterations performed at each sampling instant and establish sufficient conditions for robust constraint satisfaction.
We also provide a proof which extends the classical preconditioned fixed-point type analysis of Newton's method, see \cite[Section 5.4.2]{kelley1995iterative}, to the setting of generalized equations and establish conditions under which discrete time optimal control problems with polyhedral constraints are strongly regular. The latter property is important since it is a sufficient condition for Lipschitz continuity of the optimal value function and thus for robust stability.

This paper builds upon the results in \cite{ecc_sosspc2019} which analyzes the stability of MPC implemented using a suboptimal semismooth predictor-corrector (SSPC) method. Specifically, we generalize the previous results for suboptimal SSPC to a wide class of optimizers which are at least q-linearly convergent. Moreover, we consider external disturbances in our analysis and analyze several different variants of the RTI scheme.

The layout of the paper is as follows. We review pertinent notation and concepts in Section~\ref{ss:preliminaries} then describe the problem setting and the class of optimization algorithms we consider in Sections~\ref{ss:setting} and \ref{ss:optimization}. We establish the ISS properties of the optimization algorithms and of the coupled plant-optimizer system in Sections~\ref{ss:opt_iss} and \ref{ss:SOMPC_ISS}. Next, we discuss the strong regularity assumption in Section~\ref{ss:regularity}, we discuss relevant SQP methods in Section~\ref{ss:sequential_quadratic_programming}, and we illustrate how they fit into our optimization framework in Section~\ref{ss:time_dist_sqp}. Finally, we present simulation results in Section~\ref{ss:numerical_results}.

\section{Preliminaries} \label{ss:preliminaries}

We denote by $\ints_{+(+)}$ the non-negative (positive) integers and by $\reals_{+(+)}$ the non-negative (positive) reals. For a discrete time system
\begin{equation} \label{eq:iss_sys}
	x_{k+1} = g(x_k,u_k),
\end{equation}
given an initial state $x_0 \in \reals^n$, and an input sequence $\mathbf{u} : \nnints \to \reals^m$ we denote its solution by $x(k,x_0,\bf{u})$. For a vector, $||\cdot||$ denotes the usual Euclidean norm, for $\mathbf{u}:\nnints \to \reals^m$ we let $||\mathbf{u}|| = \mathrm{sup}\{||u_k||: k\in \nnints \}$. We use $\overline{\lim}$ as shorthand for $\limsup$. Recall that a function $\gamma :\reals_+ \to \reals_+$ is said to be of class $\K$ if it is continuous, strictly increasing and $\gamma(0) = 0$. If it is also unbounded, then $\gamma \in \Kinf$. A function $\beta:\reals_+ \times \reals_+ \to \reals_+$ is said to be of class $\KL$ if $\beta(\cdot,s) \in \K$ for each fixed $s\geq 0$ and $\beta(r,s) \to 0$ as $s \to \infty$ for fixed $r \geq 0$. If $\gamma_1,\gamma_2 : \reals \to \reals,$ we denote their composition by $\gamma_1 \circ\gamma_2$. We use $I$ to denote the identity matrix, and use $\id:\reals \to \reals$ to denote the identity function. We denote the domain of a set-valued mapping $F$ by $\dom F$. If $A$ is a matrix then $A_i$ is its $i$th row. If $\mc{I}$ is an index set, $|\mc{I}|$ is its cardinality and $A_\mc{I}$ denotes the row wise concatenation of $A_i,~ \forall i \in \mc{I}$. For two vectors, $(a,b)$ denotes vertical concatenation. We denote the unit ball centered at $x$ by $\mc{B}(x)$, it is understood that $\mc{B} = \mc{B}(0)$. If $X$ is a closed neighbourhood of the origin, we denote its radius by $\rad{X}$, i.e., the largest $r > 0$ such that $\{x~|~ ||x|| \leq r\} \subseteq X$. The normal cone mapping of a closed-convex set $C$ is defined as
\begin{equation*}
  \mc{N}_C(v) = \begin{cases}
  \{y~ |~ y^T(w-v) \leq 0, \forall w \in C\}, & \text{if}~v \in C,\\
  \emptyset & \text{else},
  \end{cases}
\end{equation*}
and set addition/subtraction is defined as
\begin{equation*}
	A \pm B = \{y~|~y = a\pm b,~ a\in A, b\in B\}.
\end{equation*}
We make extensive use of the concept of input-to-state stability\cite{jiang2001input}. Since MPC is a constrained control technique, meaning that it is intrinsically not ``global'', it is natural to consider a local variant.
\begin{dfn}[LISS\cite{jiang2004nonlinear}] \label{def:LISS}
A system \eqref{eq:iss_sys} is said to be Locally Input-to-State Stable (LISS) if there exists $\eps> 0$, $\beta \in \KL$, and $\gamma \in \K$ such that, $\forall k \in \mathbb{Z}_+$,
\begin{equation} \label{eq:ISS_bound_max}
	||x(k,x_0,\mathbf{u})|| \leq \mathrm{max}\{\beta(||x_0||,k),\gamma(||\mathbf{u}||)\},
\end{equation}
provided $||x_0||\leq \eps$ and $||\mathbf{u}|| \leq \eps$.
\end{dfn}
\begin{dfn}[Asymptotic gain\cite{jiang2001input}]
Consider system \eqref{eq:iss_sys}, we say that it has an asymptotic gain if there exists some $\gamma \in \K$ such that
\begin{equation}
	\limsupk ||x(k,x_0,\mathbf{u})|| \leq \gamma \left( \limsupk ||u_k||\right),
\end{equation}
for all $x_0\in \reals^{n_x}$.
\end{dfn}

\section{Problem Setting and Control Strategy} \label{ss:setting}
Consider the following discrete time system,
\begin{equation} \label{eq:ol_sys}
x_{k+1} = f_d(x_k,u_k,d_k),
\end{equation}
where $x_k\in \mc{X} \subset \reals^{n_x}$, $u_k\in \mc{U} \subset \reals^{n_u}$ and $d_k \in \mc{D} \subset \reals^{n_d}$ denote the state, input, and disturbance. Throughout this paper we assume full state feedback and that the following holds.
\begin{ass}\label{ass:sys_properties}
The function $f_d$ in \eqref{eq:ol_sys} is twice continuously differentiable in its first two arguments, Lipschitz continuous in the third, and $f_d(0,0,0)=0$. Moreover, the sets $\mc{X}$, $\mc{U}$ and $\mc{D}$ are compact and contain the origin.
\end{ass}

We wish to control \eqref{eq:ol_sys} using MPC and thus consider an OCP of the following form,
\begin{subequations} \label{eq:OCP}
\begin{align}
\underset{\xi,\mu}{\mathrm{min.}}~~& \phi(\xi,\mu) = V_f(\xi_N) + \sum_{i = 0}^{N-1} l(\xi_i,\mu_i),\\
\mathrm{s.t.}~~ &\xi_{i+1} = f_d(\xi_i,\mu_i,0), \quad i= 0,\ldots, N-1,\\
& \xi_0 = x,~~ \xi_N\in \mc{X}_f,\label{eq:OCP_param}\\
& (\xi_i,\mu_i) \in \mc{Z}, \quad \quad \quad \quad i = 0,\ldots, N-1.
\end{align}
\end{subequations}
where $N \in \ints_{++}$ is the horizon length, $\mc{Z} \subseteq \mc{X}\times\mc{U}$ are the constraints, $\mc{X}_f \subseteq \mc{X}$ is the terminal state constraint, $\xi = (\xi_0,\xi_1,\dots,\xi_N)$ is the state sequence, and $\mu = (\mu_0,\mu_1,\dots,\mu_{N-1})$ is the control sequence. The OCP \eqref{eq:OCP} is parameterized by the measured system state $x$. We impose the following conditions on \eqref{eq:OCP} to ensure that it is well posed and can be used to construct a stabilizing control law for \eqref{eq:ol_sys}.
\begin{ass} \label{ass:ocp_continuity}
All functions in \eqref{eq:OCP} are twice continuously differentiable in their arguments and their second derivatives are Lipschitz continuous.
\end{ass}

\begin{ass}\label{ass:terminal_trio}
The stage cost satisfies $l(0,0) = 0$, and there exists $\alpha_l \in \Kinf$ such that $\alpha_l(||x||) \leq l(x,u)$ for all $(x,u)\in \mc{Z}$. The terminal set $\mc{X}_f$ is a subset of $\mc{X}$, contains the origin in its interior, and is an admissible control invariant set for \eqref{eq:ol_sys}, i.e., for all $x\in \mc{X}_f$ there exists $u$ such that $f_d(x,u,0)\in \mc{X}_f$ and $(x,u) \in \mc{Z}$. $V_f$ is a Control Lyapunov Function for \eqref{eq:ol_sys} with $d = 0$, such that, 
\begin{equation*}
	\underset{u}{\mathrm{min}} \{V_f(x^+) - V_f(x) + l(x,u)~ | ~ (x,u) \in \mc{Z}, x^+ \in \mc{X}_f\} \leq 0,
\end{equation*}
for all $x \in \mc{X}_f$, where $x^+ = f_d(x,u,0)$.
\end{ass}


Denote by
\begin{equation}
\Gamma = \{x\in \mc{X}~|~ \text{\eqref{eq:OCP} is feasible}\},
\end{equation}
the set of feasible parameters. Under Assumptions \ref{ass:sys_properties} and \ref{ass:ocp_continuity}, the set $\Gamma$ is compact and, for all $x\in \Gamma$, a minimum of \eqref{eq:OCP} exists \cite{mayne2000constrained}. The ideal/optimal MPC feedback policy is then
\begin{equation} \label{eq:mpc_feedback}
	\kappa(x) = \mu_0^*(x),
\end{equation}
where $\mu^*(x)$ is a global minimizer of \eqref{eq:OCP}. To address the effects of incomplete optimization, we consider the perturbed closed loop system 
\begin{equation} \label{eq:cl_dyn}
x_{k+1} = f(x_k,\Delta u_k,d_k) \coloneqq f_d(x_k,\kappa(x_k)+\Delta u_k,d_k),
\end{equation}
where the control signal is corrupted by an additive disturbance that represents suboptimality, i.e., $u_k = \kappa(x_k) + \Delta u_k$. Before stating the stability properties of \eqref{eq:cl_dyn}, recall the following notion.
\begin{dfn}[RPI set \cite{limon2009input}] Given suitable sets $\Delta \mc{U}\subseteq \reals^{n_u}$ and $\mc{D}\subseteq \reals^{n_d}$, a set $\Omega \subseteq \reals^{n_x}$ is a Robust Positively Invariant (RPI) set for system \eqref{eq:cl_dyn} if $f(x,\Delta u,d) \in \Omega$ for all $x\in \Omega$, $\Delta u \in \Delta \mc{U}$, $d \in \mc{D}$. In addition, if $\Omega \subseteq \{x~|~(x,\kappa(x))\in \mc{X}\times \mc{U}\}$, then $\Omega$ is called an admissible RPI set.
\end{dfn} 

The following theorem summarizes the LISS properties of nominal MPC.
\begin{thm} \label{thm:MPC_stab} \cite[Theorem 4]{limon2009input}
Let Assumptions~\ref{ass:sys_properties} - \ref{ass:terminal_trio} hold and suppose that $\phi^*(x)$, the optimal value function for \eqref{eq:OCP}, is uniformly continuous. Then, the closed-loop system~\eqref{eq:cl_dyn} is LISS with respect to $(\Delta u,d)$ on a non-empty RPI set $\Omega \subseteq \Gamma$. Moreover, there exist $c_1,c_2 > 0$ such that if $\rad{\Delta \mc{U}} \leq c_1$ and $\rad{\mc{D}} \leq c_2$, then $\Omega$ is an admissible RPI set for \eqref{eq:cl_dyn}.
\end{thm}

In this paper we consider the situation where not enough computational resources are available to accurately solve \eqref{eq:OCP} at each sampling instant. Instead we will approximately track solution trajectories of \eqref{eq:OCP} as the measured state $x$ in \eqref{eq:OCP_param} varies over time. To track the solution trajectories, we use an appropriate iterative optimization algorithm, e.g. SQP, which is \textit{warmstarted} at time instance $t_k$ with the approximate solution from $t_{k-1}$. In this way we construct a dynamic system,
\begin{equation} \label{eq:op_dyn}
	z_{k} = \mc{T}_\ell(z_{k-1},x_k),
\end{equation}
where $z$ is an estimate of the primal-dual solution of \eqref{eq:OCP} and $\mc{T}_\ell$ represents a fixed number of optimizer iterations ($\mc{T}_\ell$ is formally defined in \eqref{eq:Tl_def}). This leads to the interconnected system illustrated in Figure~\ref{fig:subopt_mpc_fig}. The objective of this paper is to analyze the interconnection between the plant and the dynamic controller from a systems theoretic perspective.

\begin{rmk}  \label{rmk:global_opt}
This paper focuses on a common nominal MPC formulation for the sake of clarity. However, our analysis is applicable to any MPC formulation for which it is possible to prove LISS e.g., formulations that employ exact penalty functions \cite{di1989exact} or robust MPC formulations \cite{limon2009input}. Moreover, note that in \eqref{eq:mpc_feedback} the MPC feedback law is defined using a global optimum of \eqref{eq:OCP}. This requirement is not intrinsic to our analysis but rather an artifact of the specific MPC formulation. Our analysis is performed relative to a nominal ``ideal'' MPC feedback law. If the nominal feedback law is input-to-state stable our analysis is applicable regardless of whether the nominal feedback law is globally optimal or not. For example, one could use the dual-mode MPC formulation in \cite{scokaert1999suboptimal} which does not require global optimality.

\end{rmk}

\section{Optimization Strategy} \label{ss:optimization}
In this section we describe the class of optimization algorithms considered in this paper. We start in an abstract setting to clarify which properties are essential to our analysis. Later in Sections~\ref{ss:sequential_quadratic_programming} and \ref{ss:time_dist_sqp} we will illustrate how SQP fits into this framework. 

Suppose the first order necessary conditions for \eqref{eq:OCP} can be written as a parameterized Generalized Equation (GE) of the form
\begin{equation} \label{eq:GE}
	F(z,x) + \mc{N}_K(z) \ni 0,
\end{equation}
where $\mc{N}_K: \reals^n \rightrightarrows \reals^n$ is the normal cone mapping\footnote{See \cite{dontchev2009implicit} for more background on set-valued and normal cone mappings.} of a closed, convex set $K \subseteq \reals^n$, $F:\reals^n \times \Gamma \to \reals^n$ is a function, $z\in \reals^n$ are the optimization variables and $x\in \Gamma$ is the parameter. Its solution mapping is
\begin{equation} \label{eq:solution_mapping}
	S(x) = \{z ~|~ F(z,x) + \mc{N}_K(z) \ni 0\},
\end{equation}
which can be set valued. Because \eqref{eq:GE} are necessary conditions for \eqref{eq:OCP} and a minimizer of \eqref{eq:OCP} exists under Assumptions \ref{ass:sys_properties} and \ref{ass:ocp_continuity} \cite{mayne2000constrained} we have that $\dom S = \Gamma$.

Many optimization algorithms are designed to ``solve'' necessary conditions. We thus associate \eqref{eq:GE} with an iterative optimization algorithm of the form 
\begin{equation} \label{eq:opt_def}
	z_{i+1} = \mc{T}(z_i,x,i),
\end{equation}
where $\mc{T}:\reals^n \times \Gamma \times \nnints \to \reals^n$. Multiple iterations of the algorithm can be represented by the action of the function $\mc{T}_\ell(z,x): \reals^n \times \Gamma \to \reals^n$ which is defined recursively by the sequence
\begin{equation} \label{eq:Tl_def}
	\mc{T}_{\ell}(z,x) = \mc{T}(\mc{T}_{\ell-1}(z,x,\ell-1),x,\ell),
\end{equation}
where $\ell \in \ints_{++}$ is the iteration number and $\mc{T}_0(z,x) = z$.

\begin{rmk} 
The optimality conditions of \eqref{eq:OCP} can be written in multiple forms depending on the choice of \eqref{eq:GE} and \eqref{eq:opt_def}. For example, $\xi$ can either be treated as a decision variable or a function of $\mu$. As a result, the definition of $z$ is not unique and is chosen by the designer of the optimization algorithm. Specifically, the vector $z$ always includes the control sequence, but may also include the state sequence and/or the Lagrange multipliers associated with equality (dynamic) or inequality constraints.
\end{rmk}

In this paper we consider algorithms that converge at least q-linearly for a fixed parameter $x$. The following definition formalizes this notion.
\begin{dfn}[At least q-linear convergence] \label{dfn:qlinear_convergence}
For any $x\in \Gamma$ and $z^* \in S(x)$ an optimization algorithm $\mc{T}$ converges to $z^*$ if there exists $\veps > 0$ such that
\begin{equation}
	\lim_{\ell \to \infty} \mc{T}_\ell(z,x) = z^* 
\end{equation}
for all $z \in \veps \mc{B}(z^*)$. If there exists $\eta > 0$ and $q \geq 1$ such that
\begin{equation} \label{eq:one_step_error_bound}
	||\mc{T}_\ell(z,x) - z^*|| \leq \eta ||\mc{T}_{\ell-1}(z,x) - z^*||^q,
\end{equation}
for all $\ell > 0$ and $\eta \veps^{q-1} < 1$ then $\mc{T}$ is said to converge at least q-linearly over $\Gamma$.
\end{dfn}
\begin{rmk}
The necessary conditions \eqref{eq:GE} may be satisfied at stationary points or local maxima as well as local minima. We perform a local analysis that allows us to exclude those points.
\end{rmk}

Next we consider the regularity properties of \eqref{eq:solution_mapping}, we will use of the following regularity condition for generalized equations.
\begin{dfn}[Strong Regularity\cite{robinson1980strongly}]
A set-valued mapping $\Psi:\reals^n \rightrightarrows \reals^n$ is said to be strongly regular at $x$ for $y$ if $y \in \Psi(x)$ and there exist neighborhoods $U$ of $x$ and $V$ of $y$ such that the truncated inverse mapping $\tilde{F}^{-1}: V\mapsto F^{-1}(V)\cap U$ is single-valued, i.e., a function, and is Lipschitz continuous on $V$.
\end{dfn}
Strong regularity reduces to non-singularity of the Jacobian matrix if $\Psi$ is a continuously differentiable function. Our main regularity assumption follows. It establishes that any solution trajectories are Lipschitz continuous.  In Section~\ref{ss:regularity}, we discuss conditions for strong regularity for specific instances of \eqref{eq:GE}.
\begin{ass} \label{ass:pointwise_strong_reg}
All points $(z,x)$ satisfying $z\in S(x)$ that correspond to minimizers are strongly regular $\forall x \in \Gamma$.
\end{ass}

The following theorem shows that Assumption~\ref{ass:pointwise_strong_reg} ensures that the notion of tracking a local solution trajectory is well defined. 
\begin{thm} \label{thm:sol_traj}\cite{dontchev2013euler} Let the parameter $x \in \Gamma$ be a Lipschitz continuous function of $t \geq 0$. Then each solution trajectory $z(t) \in S(x(t))$ is isolated and Lipschitz continuous.
\end{thm}
A local optimization algorithm, such as SQP, can be used to track a specific ``isolated branch'' of the solution mapping. The branch is implicitly selected through the choice of initial guess supplied to the algorithm. Some MPC formulations require global optima while some do not, as discussed in Remark~\ref{rmk:global_opt}. In practice, local methods like SQP are often used regardless due to the prohibitive computational complexity of global methods.
\begin{rmk} \label{rmk:fit_the_framework}
To summarize, an algorithm/optimality condition pair fits in our framework if:
\begin{itemize}
	\item The optimality conditions can be written in the form \eqref{eq:GE} and satisfy Assumption~\ref{ass:pointwise_strong_reg}.
	\item The algorithm can be written in the form \eqref{eq:opt_def}.
	\item The algorithm is at least q-linearly convergent.
\end{itemize}
\end{rmk}
\section{LISS of time-distributed optimization} \label{ss:opt_iss}
Consider the application of TDO to problem \eqref{eq:OCP}. In a real-time setting it is only possible to perform a finite number of iterations per sampling instant which we denote by $\ell$. In this situation the optimizer can be viewed as a dynamical system of the form,
\begin{subequations}  \label{eq:optimizer_sys}
\begin{gather}
z_k = \mc{T}_\ell(z_{k-1},x_k),\\
u_k = \Xi z_k,
\end{gather}
\end{subequations}
where $\Xi$ is a surjective matrix that selects $\mu_0$ from the solution estimate, i.e., $\kappa(x) = \Xi \bar{s}(x)$ where
\begin{equation}
	\bar{s}(x) \in S(x),
\end{equation}
is an isolated single valued restriction\footnote{Our analysis is performed relative to the ideal feedback law $\kappa(x) = \Xi \bar{s}(x)$. Any choice of the restriction $\bar{s}$ that renders the origin of the closed loop system \eqref{eq:cl_dyn} LISS is admissible.} of $S$ (this is possible due to Theorem~\ref{thm:sol_traj}).

In this section we establish conditions under which \eqref{eq:optimizer_sys} is LISS. We consider the associated error system,
\begin{subequations} \label{eq:opt_err_sys}
\begin{gather} 
e_{k+1} = \mc{G}_\ell(e_k,x_k,\Delta x_k),\\
\Delta u_k = \Xi e_k,
\end{gather}
\end{subequations}
where $e_k = z_k - \bar{s}(x_k)$, and $\Delta x_k = x_{k+1} - x_k$. The error system can be explicitly constructed as follows
\begin{equation*}
	\mc{G}_\ell(e_k,x_k,\Delta x_k) = \mc{T}_\ell(e_k + \bar{s}(x_k), x_k + \Delta x_k) - \bar{s}(x_k + \Delta x_k).
\end{equation*}

\begin{lmm} \label{lmm:SL_error_bound}
Consider \eqref{eq:optimizer_sys} and its error system \eqref{eq:opt_err_sys} and suppose that $\mc{T}$ is at least q-linearly convergent. Further, let Assumption~\ref{ass:pointwise_strong_reg} hold. Then, there exists $a,\theta:\ints_{++} \to \reals_{++}$, such that the error satisfies
\begin{equation} \label{eq:SL_err_bound}
	||e_{k+1}|| \leq a(\ell) ||e_k|| +  \theta(\ell) ||\Delta x_k||,
\end{equation}
subject to the restriction
\begin{equation}\label{eq:restriction}
	||e_k|| + b ||\Delta x_k|| \leq \veps,
\end{equation}
where $\veps$ is the convergence radius in Definition~\ref{dfn:qlinear_convergence} and $b$ is the Lipschitz constant of $\bar{s}$ over $\Gamma$. Further, $a$ and $\theta$ are monotonically decreasing with $\ell$, $\lim_{\ell \to \infty} a(\ell) = 0$, $\lim_{\ell \to \infty} \theta(\ell) = 0$, and $a(\ell) \in (0,1)$.
\end{lmm}

\begin{pf*}{Proof.}
If $z_{k+1} = \mc{T}_\ell(z_k,x)$ for some fixed $x$ then the error bound \eqref{eq:one_step_error_bound} implies that
\begin{equation} \label{eq:static_SL_error_bound}
	||z_{k+1} - \bar{s}(x)|| \leq \eta^{\alpha(\ell)}||z_k - \bar{s}(x)||^{q^\ell},
\end{equation}
for all $z_k \in \veps \mc{B}(\bar{s}(x))$, d where $\alpha(\ell) = \sum_{i=0}^{\ell -1} q^i$. Now consider any $x_{k+1},x_k \in \Gamma$ and let $z_{k+1} = \mc{T}_\ell(z_k,x_{k+1})$. Then, applying \eqref{eq:static_SL_error_bound} with $x = x_{k+1}$, we obtain that
\begin{align*}
||z_{k+1} -& \bar{s}(x_{k+1})|| \leq \eta^{\alpha(\ell)} ||z_k - \bar{s}(x_{k+1})||^{q^\ell},\\
||e_{k+1}|| & \leq \eta^{\alpha(\ell)} ||z_k - \bar{s}(x_{k+1})||^{q^\ell},\\
& \leq \eta^{\alpha(\ell)} ||[z_k - \bar{s}(x_k)] - [\bar{s}(x_{k+1}) - \bar{s}(x_k)]||^{q^\ell},\\
& \leq \eta^{\alpha(\ell)}(||e_k|| + b||\Delta x_k||)^{q^\ell},
\end{align*}
where we have used that $\bar{s}$ is Lipschitz on $\Gamma$ with constant $b$ by Assumption~\ref{ass:pointwise_strong_reg}. Recall that $\veps$ denotes the convergence radius of $\mc{T}$ in $\Gamma$. A sufficient condition for the restriction $z_k \in \veps \mc{B}(\bar{s}(x_{k+1}))$ is then
\begin{align*}
&||e_k|| + b||\Delta x_k|| \leq \veps, \\
&\implies ||[z_k - \bar{s}(x_k)]|| +||[\bar{s}(x_{k+1}) - \bar{s}(x_k)]|| \leq \veps,\\
&\implies ||[z_k - \bar{s}(x_k)] - [\bar{s}(x_{k+1}) - \bar{s}(x_k)]|| \leq \veps, \\
&\implies ||z_k - \bar{s}(x_{k+1})|| \leq \veps.
\end{align*}
Continuing and imposing $||e_k|| + b||\Delta x_k|| \leq \veps$,
\begin{equation}
||e_{k+1}|| \leq\eta^{\alpha(\ell)} \veps^{q^\ell-1}(||e_k|| + b||\Delta x_k||).
\end{equation}
If $q = 1$, then $\eta^{\alpha(\ell)} \veps^{q^\ell-1} = \eta^\ell$ since $\alpha(\ell) = \sum_{i=0}^{\ell -1} q^i = \ell$ and $\veps^{q^\ell-1} = \veps^0 = 1$. Otherwise, note that
\begin{align}
\eta^{\alpha(\ell)} \veps^{q^\ell-1} &= \eta^{\frac{q^\ell-1}{q-1}} \veps^{q^\ell-1} = (\eta \veps^{q-1})^{\frac{q^\ell-1}{q-1}},
\end{align}
where we used that, for $q > 1$,
\begin{equation}
	\alpha(\ell) = \sum_{i=0}^{\ell-1} q^i = \frac{q^{\ell} -1}{q-1}.
\end{equation}
Thus,
\begin{align}
||e_{k+1}|| &\leq\eta^{\alpha(\ell)} \veps^{q^\ell-1}(||e_k|| + b||\Delta x_k||),\\
& \leq a(\ell)||e_k|| + \theta(\ell) ||\Delta x_k||,
\end{align}
where 
\begin{equation}
	a(\ell) = \begin{cases}
  \eta^\ell & \text{if}~q = 1,\\
  (\eta \veps^{q-1})^{\frac{q^\ell-1}{q-1}} & \text{if} ~ q> 1,
  \end{cases}
\end{equation}
and $\theta(\ell) = b~a(\ell)$. Since $q \geq 1$ and $\eta \veps^{q-1} < 1$ by assumption, $a(\ell) \in (0,1)$, the functions $a$ and $\theta$ are monotonically decreasing and $a(\ell),\theta(\ell) \to 0$ as $\ell \to \infty$. \qed
\end{pf*}

The following Theorem establishes the LISS properties of the error system.
\begin{thm} \label{thm:SL_iss}
Consider \eqref{eq:optimizer_sys} and its error system \eqref{eq:opt_err_sys} and suppose that $\mc{T}$ is at least q-linearly convergent. Further, let Assumption~\ref{ass:pointwise_strong_reg} hold. Then, there exists $\tau :\ints_{++} \to \reals_{++}$ such that the system is LISS if $||e_0|| \leq 0.5\veps$ and $||\Delta \mathbf{x}|| \leq \tau(\ell) \veps$, where $\tau(\ell) = 0.5~(\sigma(\ell) + b)^{-1}$ and $\veps, b$, and $\sigma$ are defined in Lemma~\ref{lmm:SL_error_bound}. Further, the asymptotic gain of \eqref{eq:opt_err_sys} is of the form
\begin{equation}
	\gamma_\ell(s) = 2\sigma(\ell)s_1 + 0\cdot s_2
\end{equation}
where $s_1$ and $s_2$ correspond to the $\Delta x$ and $x$ inputs, respectively, and $\sigma(\ell) \to 0 $ monotonically as $\ell\to\infty$.
\end{thm}
\begin{pf*}{Proof.}
Given Lemma~\ref{lmm:SL_error_bound}, if \eqref{eq:restriction} holds for all time instants leading up to $k - 1$, it follows by 
direct computation (see e.g. \cite[Example 3.4]{jiang2001input}) that
\begin{subequations} \label{eq:SL_iss1}
\begin{align}
||e_k|| &\leq a(\ell)^k ||e_0|| + \theta(\ell)\sum_{j=0}^{k} a(\ell)^{k-j}||\Delta x_j||,\\
& \leq a(\ell)^k ||e_0|| + \sigma(\ell) ||\Delta \mathbf{x}||,
\end{align}
\end{subequations}
where $\sigma(\ell)= \theta(\ell)/(1-a(\ell))$. To ensure that \eqref{eq:restriction} holds, we first consider the case $k=0$ and note that\footnote{Recall that $a + b \leq \max(2a,2b)$ for any two scalars.}
\begin{align}
	||e_0|| + b||\Delta x_0|| &\leq \max\{2||e_0||, 2 b ||\Delta \mathbf{x}|| \} ,\\
	&\leq \max\{\veps, \veps\cdot b/ (\sigma(\ell) + b)\} = \veps.
\end{align}
Next, assuming \eqref{eq:restriction} holds for $k-1$ and recalling that $a(\ell) < 1$, we can apply \eqref{eq:SL_iss1} at iteration $k$ to show that
\begin{align*}
||e_k|| + b||\Delta x_k|| & \leq a(\ell)^k ||e_0|| + \sigma(\ell) ||\Delta \mathbf{x}|| + b||\Delta x_k||,\\
& \leq a(\ell)^k ||e_0|| + (\sigma(\ell)+b)||\Delta \mathbf{x}||,\\
& \leq \max\{2 a(\ell)^k ||e_0||, 2 (\sigma(\ell)+b)||\Delta \mathbf{x}||\},\\
& \leq \max\{\veps, \veps\} = \veps,
\end{align*}
thus ensuring that \eqref{eq:restriction} also holds at $k$ due to the restriction $||e_0|| \leq 0.5 \veps$ and $||\Delta \mathbf{x}|| \leq \tau(\ell)\veps$. Since \eqref{eq:SL_iss1} recursively enforces its restrictions, we obtain
\begin{equation} \label{eq:SLiss2}
	||e_k|| \leq \max\{2 a(\ell)^k ||e_0||, 2\sigma(\ell)||\Delta \mathbf{x}||\},
\end{equation}
which directly establishes LISS with $\beta_\ell(s,k) = 2 a(\ell)^k s$ and $\gamma_\ell(||(\Delta \mathbf{x}, \mathbf{x})||) = 2\sigma(\ell) ||\Delta \mathbf{x}|| + 0 \cdot ||\mathbf{x}||$. The remaining claims follow from the expression $\sigma(\ell) = \theta(\ell)/(1-a(\ell))$ since $a(\ell), \theta(\ell) \to 0$ monotonically as $\ell \to \infty$.\qed
\end{pf*}

\section{LISS properties of suboptimal MPC} \label{ss:SOMPC_ISS}
\begin{figure}[ht]
	\centering
	\includegraphics[width=0.95\columnwidth]{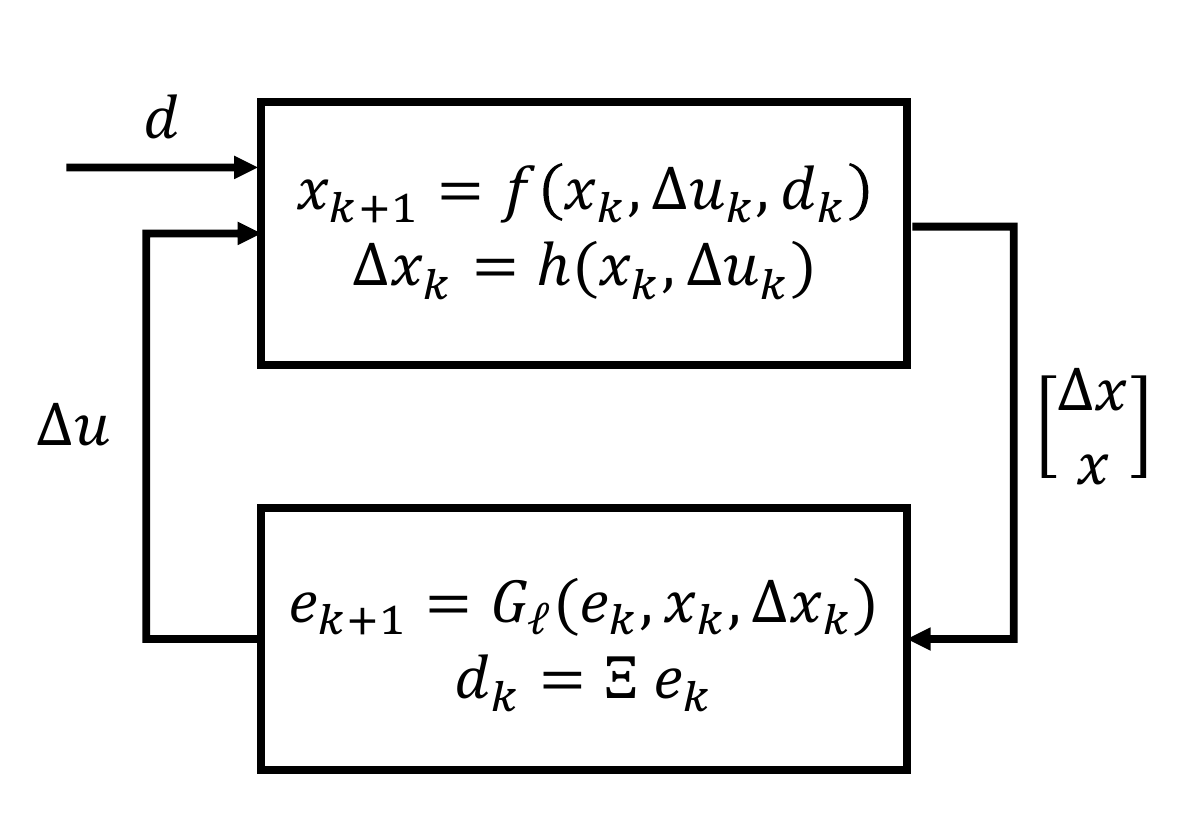}
	\caption{The coupled plant-optimizer error system.}
	\label{fig:iss_diagram}
\end{figure}

Theorem~\ref{thm:SL_iss} establishes sufficient conditions under which an at least q-linearly convergent optimizer, viewed as a dynamic system, is LISS. It also establishes that the asymptotic gain of \eqref{eq:opt_err_sys} can be made arbitrarily small by increasing the number of iterations. Since the closed-loop system \eqref{eq:cl_dyn} is itself LISS, we can derive sufficient conditions under which the coupled system, as shown in Figure~\ref{fig:iss_diagram}, is LISS with respect to the disturbance input $d$. 

\begin{thm} \label{thm:iss_proof}
Consider the dynamical systems
\begin{subequations}  \label{eq:plant_opt}
\begin{align} \label{eq:sys1}
	\Sigma_1: &\begin{cases}
	~~ x_{k+1} = f(x_k,\Delta u_k,d_k),\\
	~~ \Delta x_k = h(x_k,\Delta u_k,d_k),
	\end{cases}\\
	\Sigma_2:&\begin{cases}  \label{eq:sys2}
	~~e_{k+1} = \mc{G}_\ell(e_k,x_k,\Delta x_k),\\
	~~ \Delta u_k = \Xi e_k
	\end{cases}
\end{align}
\end{subequations}
where $h(x,\Delta u,d) = f(x,\Delta u,d) - x$ and $f$ is defined in \eqref{eq:cl_dyn}. Let the optimization algorithm used to construct $\mc{G}_\ell$ satisfy \eqref{eq:one_step_error_bound} and let Assumptions~\ref{ass:sys_properties} - \ref{ass:pointwise_strong_reg} hold. Then, there exists $\ell^* > 0$ such that if $\ell \geq \ell^*$ the interconnected system \eqref{eq:plant_opt} is LISS with respect to the input $d$.
\end{thm}
\begin{pf*}{Proof.}
Under Assumptions~\ref{ass:sys_properties} - \ref{ass:pointwise_strong_reg}, Theorem \ref{thm:MPC_stab} holds\footnote{Recall that Assumption~\ref{ass:pointwise_strong_reg} is sufficient for Lipschitz continuity of the optimal value function $\phi^*(x)$.}. Thus system $\Sigma_1$ is LISS, meaning that there exist asymptotic gains $\gamma_1,\gamma_2 \in \K$ such that
\begin{multline} \label{eq:sys1_gain}
	\limsupk ||x_k|| \leq \gamma_1\left(\limsupk ||\Delta u_k||\right) + \gamma_2\left(\limsupk ||d_k||\right),
\end{multline}
for suitably restricted $d_k\in\mathcal{D}$ and $\Delta u_k\in\Delta\mathcal{U}$. Let $L$ denote the Lipschitz constant of $h$, then
\begin{equation}
	||\Delta x_k|| \leq L ||x_k|| + L ||\Delta u_k|| + L ||d_k||,
\end{equation}
combining this with \eqref{eq:sys1_gain} we obtain that
\begin{equation} \label{eq:iss2}
\limsupk ||\Delta x_k|| \leq  \gamma_3\left(\limsupk ||\Delta u_k||\right) + \gamma_4\left(\limsupk ||d_k||\right),
\end{equation}
where $\gamma_3 = L (\gamma_1  + \id)$, and $\gamma_4 = L (\gamma_2 + \id)$. Similarly, due to  Theorem~\ref{thm:SL_iss}, there exists $\veps > 0$ and positive functions $\sigma$ and $\tau$ such that
\begin{equation}
	\limsupk ||e_k|| \leq \sigma(\ell) \limsupk ||\Delta x_k||,
\end{equation}
given $||\Delta \mathbf{x}||\leq\tau(\ell)\veps$. Therefore, it follows from \eqref{eq:sys2} that
\begin{subequations} \label{eq:iss1}
\begin{align} 
	\limsupk ||\Delta u_k|| &\leq ||\Xi|| \limsupk ||e_k||,\\
	& \leq \sigma(\ell)||\Xi|| \limsupk ||\Delta x_k||.
\end{align}
\end{subequations}
Combining \eqref{eq:iss1} with \eqref{eq:iss2} we obtain that
\begin{multline}
\limsupk ||\Delta u_k|| \leq \sigma(\ell)||\Xi|| \gamma_3 \left(\limsupk ||\Delta u_k||\right) + \\\sigma(\ell)||\Xi|| \gamma_4\left(\limsupk ||d_k||\right).
\end{multline}
Thus, if the contraction property
\begin{equation} \label{eq:small_gain}
	||\Xi|| \sigma(\ell)\gamma_{3}(s) \leq s
\end{equation}
is satisfied for all $s\in[0,\rad\Delta\mathcal{U}]$, \eqref{eq:plant_opt} is LISS with suitable restrictions on the initial state and on the disturbance $d$, as detailed in \cite[Theorem 2]{teel1996}. Note that, since $\Delta u = \Xi e$, we have $\rad \Delta \mc{U} \leq \Xi \veps$ where $\veps$ is the convergence radius\footnote{If the optimizer is globally convergent then $\veps$ can be chosen arbitrarily. In that scenario it may be possible to obtain a stronger result using different analysis tools.} of the optimizer defined in Theorem~\ref{thm:SL_iss}. Since $\sigma(\ell) \to 0$ monotonically as $\ell \to \infty$, the existence of $\ell^* < \infty$ such that \eqref{eq:small_gain} is satisfied follows from the finiteness of $\gamma_{3}$, $\rad\Delta\mathcal{U}$, and $||\Xi||$.
\qed
\end{pf*}
Theorem~\ref{thm:iss_proof} establishes conditions under which the interconnected plant-optimizer system is LISS. However, this result does not provide any information about the set of admissible initial conditions and does not consider constraint satisfaction. By noting that the ideal MPC feedback law admits a robust positively invariant set, we can extend our result by deriving sufficient conditions for constraint satisfaction.
\begin{thm} \label{thm:constraint_satisfaction}
Suppose that the assumptions of Theorem~\ref{thm:iss_proof} hold so the interconnected system \eqref{eq:plant_opt} is LISS. Let $\Omega$ denote the admissible RPI set in Theorem~\ref{thm:MPC_stab}, let $\gamma_\ell(s) = 2\sigma(\ell) s$ denote the asymptotic gain of \eqref{eq:opt_err_sys}, and let $(x_k,e_k) = (x(k,x_0,\mathbf{d}),e(k,e_0,\mathbf{d}))$ denote the closed-loop trajectory of \eqref{eq:plant_opt} for some initial condition $(x_0,e_0)$ and disturbance sequence $\mathbf{d}$. Then, if the disturbances are sufficiently small, there exists $\bar{\ell} \geq \ell^*$ and $\delta > 0$ such that, if $||e_0|| \leq \delta$ and $x_0\in \Omega$, then $x_k \in \Omega$ for all $k\geq 0$.
\end{thm}
\begin{pf*}{Proof.}
Due to Theorem~\ref{thm:MPC_stab}, given a sufficiently small disturbance set $\mc{D}$, there exists a neighbourhood of the origin $\Delta \mc{U}$ such that, if $\Delta u_k \in \Delta \mc{U}, ~\forall k \geq 0$ and $x_0 \in \Omega$, then $x_k \in \Omega,~\forall k\geq 0$. Since $\Delta u = \Xi e$ for a surjective matrix $\Xi$, there exists $\rho> 0$ such that, if $||e_k|| \leq \rho$, then $\Delta u_k \in \Delta \mc{U}$. Given the restriction $||\Delta \mathbf{x}|| \leq \tau(\ell) \veps$ and $||e_0|| \leq 0.5 \veps$, where $\tau$ and $\veps$ are defined in Theorem~\ref{thm:SL_iss}, it follows from \eqref{eq:SLiss2} that $||e_k|| \leq \rho$ can be imposed by enforcing  $||e_0|| \leq 0.5\rho$ and $2\sigma(\ell) ||\Delta \mathbf{x}|| \leq \rho$. To enforce $2\sigma(\ell) ||\Delta \mathbf{x}|| \leq \rho$, we note that the set $\Omega$ is bounded \cite[Theorem 4]{limon2009input}, thus implying that $\Delta x \in \Delta \Omega = \Omega - \Omega$ (see Section~\ref{ss:preliminaries} for a definition of set subtraction) is bounded by 
\begin{equation}
	\bar{s} = \underset{w\in \Delta \Omega}{\sup}~||w|| < \infty.
\end{equation}
Since $\bar{s}$ is finite and $\sigma(\ell) \to 0$ monotonically as $\ell \to \infty$, there exists $\ell_1$ such that $2\sigma(\ell_1)\bar{s} \leq \rho$. Moreover, since $\bar{s}$ is finite and $\tau(\ell) \to \infty$ monotonically as $\ell \to \infty$ there exists $\ell_2$ such that $\bar{s} \leq \tau(\ell_2)\veps$. Thus, letting $\delta = 0.5 \max\{\rho,\veps\}$ and $\bar{\ell} = \max(\ell^*,\ell_1,\ell_2)$, it follows that the system is LISS with restrictions on the initial conditions $x_0\in\Omega$ and $\|e_0\|\leq\delta$, as well as restrictions on the external disturbance $d\in\mathcal{D}$.\qed
\end{pf*}
Theorem~\ref{thm:constraint_satisfaction} establishes that, if enough computational resources are available and the initial solution guess is sufficiently accurate, then TDO recovers the robustness properties of optimal MPC.

\begin{rmk}
The results presented in this section are quite general: as long as the MPC formulation is LISS, the solution mapping of the OCP is strongly regular, and the convergence rate of the iterative solver is at least $q$-linear, Theorems \ref{thm:iss_proof} and \ref{thm:constraint_satisfaction} prove that it is possible to achieve robust stability and constraint satisfaction by performing a limited number of solver iterations per time step. Due to the generality of the framework, however, the actual values we obtain for $\ell^*$ and $\bar\ell$ are likely to be conservative and would be ill-suited for, e.g., complexity certification as in \cite{rubagotti2014stabilizing}, which, it should be noted, only considers the linear case. Despite this drawback, our results significantly extend the existing analysis of the RTI scheme \cite{diehl2005nominal} when our framework is applied to TD-SQP. Complexity certification is significantly more challenging in the nonlinear case due to the nonconvexity of the OCPs and is left to future work. 
\end{rmk}

\section{Conditions for Strong Regularity}\label{ss:regularity}
The main results in this paper are all predicated upon
Assumption 4, that for each parameter value $x\in\Gamma$ the
solution mapping of the OCP is strongly regular. In this section we discuss some common settings and derive strong regularity conditions for each.

\subsection{Closed Convex Constraint Sets} \label{ss:convex_constraint_sets}
If the constraint sets $\mc{Z}$ and $\mc{X}_f$ in \eqref{eq:OCP} are closed and convex, it is possible to write the optimality conditions without introducing dual variables for the inequality constraints. In particular, we can express \eqref{eq:OCP} compactly as
\begin{equation} \label{eq:OCP_convex_compact}
	\underset{w \in W}{\mathrm{min}}~~ \phi(w),~~\mathrm{s.t.}~~ g(w,x) = 0,
\end{equation}
where $W  = \mc{Z} \times \mc{Z} \ldots \times \mc{X}_f$ and $w = (\xi_0,\mu_0,\ldots, \xi_N) \in \reals^p$. The Lagrangian associated with \eqref{eq:OCP_convex_compact} is
\begin{equation}
	\mc{L}(w,\lambda,x) = \phi(w) + \lambda^Tg(w,x),
\end{equation}
where $\lambda\in \reals^{l}$ are dual variables (sometimes ``co-states''). The KKT conditions of \eqref{eq:OCP_convex_compact} are
\begin{equation} \label{eq:OCP_convex_KKT}
	\nabla_z \mc{L}(z,x) + \mc{N}_{Z}(z) \ni 0,
\end{equation}
where $z = (w,\lambda)$ and $Z = W \times \reals^{l}$. Note that \eqref{eq:OCP_convex_KKT} can be reduced to \eqref{eq:GE} by choosing $F = \nabla_z \mc{L}$ and $K = Z$. Our framework requires that \eqref{eq:OCP_convex_KKT} be necessary for optimality. To ensure this, we impose the following constraint qualification \cite[Theorem 6.14]{rockafellar2009variational}
\begin{equation} \label{eq:cc}
-\nabla_w g(\bar{w},\bar{x})^T y \in \mc{N}_W(\bar{w})  \implies  y = 0,
\end{equation}
for all $(\bar{w},\bar{\lambda}) \in S(\bar{x})$. The following lemma proves that\eqref{eq:cc} holds automatically in this setting.
\begin{lmm}
The constraint qualification \eqref{eq:cc} holds at all points $(z,x) \in \reals^{p+l} \times \Gamma$.
\end{lmm}
\begin{pf*}{Proof.}
The constraint qualification is implied by surjectivity of the matrix $\nabla_w g(w,x)$. Denoting $A_i = \nabla_\xi f_d(\xi_i,\mu_i,0)$, and $B_i = \nabla_\mu f_d(\xi_i,\mu_i,0)$, the surjectivity of $\nabla_w g(w,x)$ becomes the condition that for every $\xi = (\xi_0, \dots,\xi_{N})$ the system 
\begin{equation*}
	x_0 =\xi_0, ~~\zeta_{i+1}- A_i\zeta_i - B_i\nu_i = \xi_{i+1},~~i \in \ints_{[0,N-1]},
\end{equation*}
has a solution. This condition clearly holds: pick an arbitrary sequence $(\nu_0, \dots, \nu_{N-1})$ and determine $(\zeta_0,\ldots,\zeta_N)$ recursively.\qed
\end{pf*}

Before stating the conditions for strong regularity, we recall the following second order condition (which can be monitored numerically, see e.g., \cite[Section 16.2]{nocedal2006numerical}).
\begin{dfn}[SOSC] The Second Order Sufficient Condition (SOSC) is said to hold at $\bar{z} = (\bar{w},\bar{\lambda}) \in S(\bar{x})$ if
\begin{equation} \label{eq:convex_SOSC}
y^T \nabla_w^2\mc{L}(\bar{z},\bar{x}) y > 0,~~ \forall y \text{ s.t. } \nabla_w g(\bar{w},\bar{x}) y = 0.
\end{equation}
\end{dfn}

\subsubsection{Convex Control Constraints} \label{ss:convex_control_constraints}
If only convex control constraints are present, Theorem~\ref{thm:ctrl_sr} provides sufficient conditions for strong regularity. 
\begin{thm} \label{thm:ctrl_sr} \cite[Theorem 1.2]{dontchev2019lipschitz}
Suppose that $\mc{Z} = \reals^{n_x} \times \mc{U}$, where $\mc{U}$ is closed and convex, and consider any $\bar{z} \in S(\bar{x})$. If \eqref{eq:convex_SOSC} holds, then $S$ is strongly regular at $(\bar{z},\bar{x})$. 
\end{thm}
As a result of Theorem~\ref{thm:ctrl_sr}, Assumption~\ref{ass:pointwise_strong_reg} reduces to the assumption that \eqref{eq:convex_SOSC} holds at all minimizers in $\Gamma$. In this scenario, any terminal set constraints would have to be enforced through penalty functions. 

\subsubsection{Polyhedral State and Control Constraints}
If the state and control constraints are convex polyhedra, the following theorem applies. The result was previously asserted without proof in \cite[Section 3.2]{dinh2012adjoint}, we provide a proof for completeness.
\begin{thm} \label{thm:sr_poly} Suppose that $W$ in \eqref{eq:OCP_convex_compact} is polyhedral with a representation $W = \{w~|~Mw\leq h\}$. Now consider a KKT point $\bar{z} = (\bar{w},\bar{\lambda}) \in S(\bar{x})$. If \eqref{eq:convex_SOSC} holds, then $S$ is strongly regular at $(\bar{z},\bar{x})$. 
\end{thm}
\begin{pf*}{Proof.}
Strong regularity of the nonlinear GE \eqref{eq:OCP_convex_KKT} at $(\bar{z},\bar{x})$ follows from strong regularity of its partial linearization \cite{robinson1980strongly}. This can be written as
\begin{equation} \label{eq:VIlin}
	\begin{bmatrix}
		R & G^T \\
		- G & 0
	\end{bmatrix} \begin{bmatrix}
		w \\ \lambda
	\end{bmatrix} + \begin{bmatrix}
		r \\ g
	\end{bmatrix} + \mc{N}_C(z) \ni 0,
\end{equation}
where $\hat{f} = \nabla_w\mc{L}(\bar{z},\bar{x})$, $R = \nabla_w^2\mc{L}(\bar{z},\bar{x})$, $G = \nabla_w g(\bar{z},\bar{x})$, $r = \hat{f} - R \bar{w} - G^T \bar{\lambda}$, $g = G \bar{w}$, $\Theta = \{w~|~ Gw = g, ~M w \leq h\}$ and $C = \Theta \times \reals^{l}$. Equation~\eqref{eq:VIlin} is an affine GE of the form,
\begin{equation}
	A z + a + \mc{N}_C(z) \ni 0,
\end{equation}
to which we apply \cite[Theorem 2E.6]{dontchev2009implicit} to establish strong regularity of the mapping $A+ \mc{N}_C$. This requires
\begin{equation} \label{eq:polySR_conditions}
	z \in \mc{E}^+,~Az\perp \mc{E}^-,~ z^T A z \leq 0 \Rightarrow z = 0,
\end{equation}
where $\mc{E}^+ = \mc{E}-\mc{E}$, $\mc{E}^- = \mc{E} \cap -\mc{E}$, and
\begin{equation*}
	\mc{E} =  \{(w,\lambda)~|~ Gw = 0, M_i w \leq 0~i\in \mc{A}(\bar{w}), \hat{f}^T w = 0\},
\end{equation*}
is the critical cone\footnote{See \cite[Section 2E]{dontchev2009implicit} for more details on critical cones. We've simplified the expression for $\mc{E}$ using \cite[Theorem 2E.3]{dontchev2009implicit} and \eqref{eq:OCP_convex_KKT}.} of $C$ at $\bar{z}$. Next, note that $\mc{E} \subseteq \mc{E}^+ \subset \ker G \times \reals^{l}$ thus, by the second order condition, $y^T R y = y^T \nabla_w^2 \mc{L}(\bar{z},x)y > 0$ for all $y \in \ker G$. Thus 
\begin{equation}
	z^T A z = w^T R w > 0, \forall w \in \mc{E}^+ ,
\end{equation}
which implies that $z^T A z \leq 0 \implies z = 0$ for all $z \in \mc{E}^+$. As a result, \eqref{eq:polySR_conditions} is satisfied and \eqref{eq:OCP_convex_KKT} is strongly regular. \qed
\end{pf*}
Thus, as in the case of convex input constraints, Assumption~\ref{ass:pointwise_strong_reg} reduces to the condition that \eqref{eq:convex_SOSC} holds at all minimizers in $\Gamma$. 

\subsection{Nonlinear Inequality Constraints} \label{ss:NLP}
If the constraint sets in \eqref{eq:OCP} can be expressed in the form $\mc{Z} = \{(\xi,\mu)~|~c(\xi,\mu) \leq 0\}$ and $\mc{X}_f = \{\xi~|~c_f(\xi) \leq 0\}$ for suitable twice continuously differentiable functions $c:\reals^{n_x+n_u} \to \reals^{n_c}$ and $c_f:\reals^{n_x}\to\reals^{n_{cf}}$, then \eqref{eq:OCP} can be written compactly as the following Nonlinear Program (NLP),
\begin{subequations}  \label{eq:NLP}
\begin{align}
\underset{w}{\mathrm{min.}} \quad &\phi(w),\\
\mathrm{s.t.} \quad &g(w,x) = 0,~~h(w) \leq 0,
\end{align}
\end{subequations}
where $w = (\xi,\mu) \in \reals^p$ are the decision variables. The Lagrangian associated with \eqref{eq:NLP} is
\begin{equation}
L(w,\lambda,v,x) = \phi(w) + \lambda^T g(w,x) + v^T h(w),
\end{equation}
where $\lambda \in \reals^{l}$ and $v\in \reals^{m}$ are dual variables. Its KKT conditions \cite{izmailov2014newton} are 
\begin{subequations} \label{eq:KKT}
\begin{gather}
\nabla_w L(w,\lambda,v,x) = 0,\\
-g(w,x) = 0,\\
-h(w) +\mc{N}_+(v) \ni 0, \label{eq:KKTcomp}
\end{gather}
\end{subequations}
where $\mc{N}_+$ is the normal cone mapping of the non-negative orthant. Comparing \eqref{eq:KKT} with \eqref{eq:GE} we can identify $z = (w,\lambda,v)$,
$K = \reals^{p} \times \reals^{l} \times \reals_{\geq 0}^m$, and
\begin{equation}
	F(z,x) = \begin{bmatrix}
		\nabla_w L(w,\lambda,v,x)\\
		-g(w,x)\\
		-h(w)
	\end{bmatrix}.
\end{equation}
To ensure that \eqref{eq:KKT} are necessary for optimality, as required by our framework, we need to impose a constraint qualification on \eqref{eq:NLP}.
\begin{dfn}[LICQ] \label{dfn:LICQ}
The Linear Independence Constraint Qualification (LICQ) is said to hold at $(\bar{z},\bar{x})$ if
\begin{equation*}
  \text{rank}~\begin{bmatrix}
  \nabla_w g(\bar{w},\bar{x})\\
  [\nabla_w h(\bar{w})]_i
  \end{bmatrix}
   = l + |\mc{A}(\bar{w})|,~ i \in \mc{A}(\bar{w}),
\end{equation*}
where $\mc{A}(w) = \{i\in 1~...~ q~|~ h_i(w) =0\}$ is the set of active constraint indices and $l$ is the number of equality constraints.
\end{dfn}

The following Theorem summarizes necessary and sufficient conditions for strong regularity in the context of nonlinear programming.
\begin{thm} \label{thm:strong_reg} \cite[Prop 1.27, 1.28]{izmailov2014newton}
Consider a parameterized nonlinear program of the form \eqref{eq:NLP} and let $S(x)$ be the solution mapping of its KKT conditions \eqref{eq:KKT}. A point $(\bar{z},\bar{x})$ satisfying $\bar{z} \in S(\bar{x})$ is strongly regular if it satisfies the LICQ and the strong second order sufficient condition (SSOSC)
\begin{equation*}
  y^T \nabla_w^2 L(\bar{z},\bar{x}) y> 0,~\forall y \in \mathcal{K}_+(\bar{z},\bar{x}) \setminus \{0\},
\end{equation*}
where 
\begin{multline*}
	\mathcal{K}_+(z,x) = \{y\in \reals^n~|~ \nabla_w g(\bar{w},\bar{x})y = 0,\\~ [\nabla_w h(\bar{w})]_iy = 0,~i \in \mc{A}^+(z,x)\}, 
\end{multline*}
and $\mc{A}^+(z,x) = \mc{A}(w) \cap \{i~|~ v_i > 0\}$. Moreover, if $\bar{w}$ is a local minimizer of \eqref{eq:NLP}, the LICQ and SSOSC are also necessary conditions for strong regularity. 
\end{thm}
Thus, Assumption~\ref{ass:pointwise_strong_reg} reduces to the assumption that the SSOSC and LICQ hold at all minima in $\Gamma$.

\section{Sequential Quadratic Programming}\label{ss:sequential_quadratic_programming}
Having defined under what conditions the solution mapping of the OCP is strongly regular, we investigate the convergence properties of two widely used SQP schemes to show that they can be used for TDO. To this effect, note that the OCPs \eqref{eq:OCP_convex_compact} and \eqref{eq:NLP} can both be solved using SQP. Specifically, for \eqref{eq:NLP}, given a solution estimate $z_i$, the next iterate can be computed by solving the following Quadratic Program (QP)
\begin{subequations} \label{eq:QP}
\begin{align}
\underset{\Delta w_i}{\mathrm{min.}}\quad \frac12 &\Delta w_i^T B_i \Delta w_i + \nabla_w\phi(z_{i})^T \Delta w_i,\\
\mathrm{s.t.}\quad ~~ &\nabla_w g(w_{i},x) \Delta w_i + g(w_{i},x) = 0,\\
&\nabla_w h(w_i) \Delta w_i + h(w_i) \leq 0, \label{eq:QP_ineq}
\end{align}
\end{subequations}
where $B_i$ approximates the Hessian of the Lagrangian $\nabla_w^2 L$. Specifically, if we denote the Lagrange multipliers associated with the equality and the inequality constraints by ${\pi_i}$ and ${\eta_i}$ respectively. The SQP update for \eqref{eq:NLP} is $z_{i+1} = (w_i + \Delta w_i, \pi_i,\eta_i)$. Note that \eqref{eq:QP} is fully defined by \eqref{eq:OCP_convex_compact} or \eqref{eq:NLP} except for $B_i$, which will depend on the specific SQP method. 

SQP applied to \eqref{eq:OCP_convex_compact} is similar. The SQP update becomes $z_{i+1} = (w_i + \Delta w_i, \pi_i)$, i.e., the inequality duals are removed from the iteration, $B_i$ must approximate $\nabla_z^2 \mc{L}$ instead of $\nabla_w^2 L$ and \eqref{eq:QP_ineq} becomes $M \Delta w_i + M w_i - h \leq 0$, where $M$ and $h$ satisfy $W = \{w~|~Mw \leq h\}$.
\begin{rmk}
The convex control constraints setting (Section~\ref{ss:convex_control_constraints}) technically allows for non-polyhedral convex constraints. In this case \eqref{eq:QP_ineq} would need to be replaced with a convex constraint of the form $w_i+ \Delta w_i \in W$ but otherwise no changes are necessary.
\end{rmk}  

To provide a unified formulation, we exploit that SQP can be seen as a Newton-type process for solving GEs of the form
\begin{equation} \label{eq:GE_normal_cone}
	F(z,x) + \mc{N}_K(z) \ni 0,
\end{equation}
where $z\in \reals^n$, $x\in \Gamma$, $F:\reals^n \times \Gamma \to \reals^n$ is continuously differentiable and $\mc{N}_K:\reals^n \rightrightarrows \reals^n$ is the normal cone mapping for a closed, convex set $K \subseteq \reals^n$. Newton's method applied to \eqref{eq:GE_normal_cone} is
\begin{equation} \label{eq:newton_process}
	H_i(z_{i+1} - z_i) + F(z_i,x) + \mc{N}_K(z_{i+1}) \ni 0,
\end{equation}
where the sequence $\{H_i\}$ approximates $\nabla_z F(z_i,x)$. Referring to the QP subproblem \eqref{eq:QP}, we note that
\begin{equation}
H_i = \begin{bmatrix}
B_i & \nabla_w^T g(w_i,x) & \nabla_w^T h(w_i)\\
-\nabla_w g(w_i,x) & 0 & 0 \\
-\nabla_w h(w_i) & 0 & 0
\end{bmatrix},
\end{equation}
for \eqref{eq:NLP}'' for \eqref{eq:OCP_convex_compact} simply discard the third row and column. Thus, the the sequence $\{H_i\}$ is fully determined by the Hessian approximation sequence $\{B_i\}$.

\begin{rmk}
In this paper we only consider ``undamped'' Newton methods, which are intrinsically local methods. More sophisticated implementations may include various type of regularization and/or globalization techniques such as trust regions or linesearches to enlarge the methods region of attraction. Nevertheless, undamped Newton methods are commonly used in practice, especially in the context of the RTI scheme, and the tools we develop in this paper are applicable to locally convergent algorithms. We leave the application of our tools to globalized SQP methods to future work and refer readers to e.g., \cite{nocedal2006numerical} or \cite{izmailov2014newton}, for more detailed treatments of SQP methods.
\end{rmk}

\subsection{The Josephy-Newton (JN) method} \label{ss:josephy_newton}
Using the exact Hessian of the Lagrangian results in the Josephy-Newton method. The following theorem summarizes the convergence properties of the JN method applied to \eqref{eq:GE_normal_cone}.
\begin{thm} \cite[Theorem 3.2]{izmailov2014newton}\label{thm:JN_convergence}
Let $z^* \in \bar{S}(x)$ for some fixed $x$ and suppose that Assumption~\ref{ass:ocp_continuity} holds and $(z^*,x)$ is strongly regular. Let the sequence $\{z_i\}$ be generated by repeatedly solving
\begin{equation} \label{eq:JN_process}
	\nabla_z F(z_i,x) (z_{i+1}-z_i) + F(z_i,x) + \mc{N}_K(z_{i+1}) \ni 0.
\end{equation}
Then, there exists $\bar{\eta} = \bar{\eta}(x)>0$ and $\bar{\eps} = \bar{\eps}(x) > 0$ satisfying $\bar{\eta} \bar{\eps} < 1$, such that, if $z_0 \in \bar{\eps} \mc{B}(z^*)$, then $\{z_i\}$ is unique and converges to $z^*$ q-quadratically, i.e.,
\begin{equation}
	||z_{i+1} - z^*|| \leq \bar{\eta} ||z_i - z^*||^2.
\end{equation}
\end{thm}
In general, we cannot expect $\nabla_w^2 L$ to be positive semidefinite even in the vicinity of a solution. This may make solving the QP subproblems difficult and is a well known issue in the SQP literature. A detailed discussion is outside the scope of this paper, we refer interested readers to e.g., \cite{izmailov2014newton,nocedal2006numerical,boggs1995sequential}. We will however briefly discuss two possible solutions. The first is to use an Augmented Lagrangian Hessian, i.e., to use $B = \nabla_w^2L + \rho (\nabla_w g) (\nabla_w g)^T$ for some $\rho > 0$. If the second order sufficient conditions hold, then $B$ will be convex if $\rho$ is sufficiently large \cite[Section 4.2]{izmailov2014newton}. This will shift the multipliers associated with $g$, see \cite[Section 4.2]{izmailov2014newton} for details on how to recover the original multipliers. This method may not always be numerically efficient because it can negatively impact the sparsity of $B$. The second is to use a reduced Hessian approach, see e.g., \cite{schmid1994quadratic,nocedal2006numerical}, which maintains a basis for the null space of the active constraints and solves the QPs on this reduced space. The Hessian projected onto the reduced space (the ``reduced Hessian'') is guaranteed to be positive definite in the vicinity of a solution if a second order condition holds.

\subsection{The Gauss-Newton (GN) method} \label{ss:gauss_newton}
The Gauss-Newton method is applicable when the objective function has the form $\phi(w) = ||r(w)||^2_2$ for some residual function $r$. The Hessian of the Lagrangian is then approximated by
\begin{equation}
	 B(w) = \nabla_w r(w) \nabla_w r(w)^T \approx \nabla_w^2 L(z,x) .
\end{equation}
For example if $\phi(w) = x^TQx + u^T Ru$ the GN Hessian approximation is $B = \texttt{blkdiag}(Q,R)$. The GN method has the advantage that the Hessian approximation is guaranteed to be positive semidefinite, so the QP subproblems can be solved reliably. Because of this, the GN method is widely used in practice, see e.g., \cite{houska2011acado,vukov2012experimental,gros2013economic,gros2016linear,diehl2007stabilizing,albin2018vehicle}. The GN approximation error satisfies
\begin{multline} \label{eq:GN_approx_error}
	\nabla_w^2 L(z,x) - B(w) = \\\Oh(||r(w)||) + \Oh(\sum_{i=1}^m ||\lambda_i||~||\nabla_w^2 g_i(w,x)||),
\end{multline}
so the approximation error depends on the size of the residuals and on the second derivative $g$ which is related to the nonlinearity of the dynamics. We show in Theorem~\ref{thm:GN_convergence} that it is important to approximate $\nabla_w^2 L(z^*,x)$ where $z^*\in \bar{S}(x)$. The following theorem establishes sufficient conditions for q-linear convergence of the GN method by extending the classical fixed-point type analysis of Newton's method, see \cite[Section 5.4.2]{kelley1995iterative}. The nearest analysis we found in the literature is \cite[Theorem 3.5]{dinh2012adjoint} which considers a path tracking problem rather than a fixed one.

\begin{thm} \label{thm:GN_convergence}
Fix some parameter $x \in \Gamma$, let $z^* \in S(x)$ and suppose that Assumptions~\ref{ass:ocp_continuity} and \ref{ass:pointwise_strong_reg} hold. Consider a sequence $\{z_i\}$ generated by repeatedly solving \eqref{eq:newton_process}. Further, define $e_i = z_i - z^*$ and suppose that there exist $\bar{\delta} = \bar{\delta}(x) > 0$ such that $||H_i - \nabla F(z^*,x)|| \leq \bar{\delta}$ for all $i \geq 0$. If the mapping
\begin{equation}
	J_i(z) = H_i z + \mc{N}_K(z) 
\end{equation}
is strongly regular for all $i \geq 0$, i.e., $J_i^{-1}$ is a Lipschitz continuous function with Lipschitz constant $M>0$, and $\bar{\delta} M < 1$, then there exists $\bar{\eps} = \bar{\eps}(x) > 0, $ and $L >0$ such that if $z_0 \in \bar{\eps} \mc{B}(z^*)$, then $\{z_i\}$ is unique, converges to $z^*$ q-linearly, and
\begin{equation}
	||e_{i+1}|| \leq M(\bar{\delta} + L||e_i||)||e_i|| \leq  \bar{\eta} ||e_i||,
\end{equation}
where $\bar{\eta} = \bar{\eta}(x) = M(\bar{\delta} + L\bar{\eps})$.
\end{thm}
\begin{pf*}{Proof.}
A solution, $z^* \in \bar{S}(x)$, exists for every $x \in \Gamma $ thanks to Assumption~\ref{ass:pointwise_strong_reg}; from this point forward we will suppress the dependencies on $x$ in the subsequent expressions. The GN method can be written as
\begin{equation}
	z_{i+1} = J_i^{-1} \circ G_i(z_i) = T_i(z_i),
\end{equation}
where $G_i(z) = H_i z -F(z)$; note that $z^*=T_i(z^*)$ for any choice of $\{H_i\}$. First consider
\begin{equation*}
	G_i(z_i) - G_i(z^*) = H_i(z_i - z^*) - F(z_i) + F(z^*)
\end{equation*}
\begin{multline*}
 = [\nabla F(z^*) (z_i-z^*) - F(z_i) +  F(z^*)] +\\ [(H_i - \nabla F(z^*))(z_i-z^*)].
\end{multline*}
Since $\nabla F$ is Lipschitz (Assumption~\ref{ass:ocp_continuity}) the fundamental theorem of calculus implies that there exist $L,\eps_1 > 0$ such that
\begin{equation*}
	||\nabla F(z^*) (z_i-z^*) - F(z_i) +  F(z^*)|| \leq L ||z_i - z^*||^2,
\end{equation*}
for all $z_i\in \eps_1 \mc{B}(z^*)$, so, taking norms, we obtain that
\begin{equation*}
||G_i(z_i) - G_i(z^*)|| \leq L ||z_i - z^*||^2 + \bar{\delta}||z_i-z^*||.
\end{equation*}
By assumption the mapping $J_i^{-1}$ is Lipschitz continuous so $\Delta T_i^* = ||T_i(z) - T_i(z^*)||$ satisfies
\begin{subequations} \label{eq:GN_linear_bound}
\begin{align}
	\Delta T_i^*& = ||J_i^{-1}(G_i(z_i)) - J_i^{-1}(G_i(z^*))||,\\
	&\leq M ||G_i(z_i) - G_i(z^*)||,\\
	& \leq M (\bar{\delta} + L||e_i||) ||e_i||,
\end{align}
\end{subequations}
for all $z_i \in \eps_1 \mc{B}(z^*)$. Now consider the update equation
\begin{equation*}
	||z_{i+1} - z^*|| = ||T_i(z_i) - z^*|| = ||T_i(z_i) - T_i(z^*)||,
\end{equation*}
where we have used that $z^* = T_i(z^*)$. Since $J_i$ is strongly regular, $T_i$ is a function and $\{z_i\}$ is unique. Using \eqref{eq:GN_linear_bound} we have
\begin{equation} 
	||e_{i+1}|| \leq M(\bar{\delta} + L||e_i||)||e_i||, ~~\forall e_i\in \eps_1 \mc{B}.
\end{equation}
Since $M\bar{\delta} < 1$ by assumption, it is possible to pick $\bar{\eps}\in (0,\eps_1)$ such that $\bar{\eta}  = M(\bar{\delta} + L \bar{\eps})< 1$. Then $\{z_i\}$ converges q-linearly to $z^*$ if $z_0 \in \bar{\eps}\mc{B}(z^*)$, i.e.,
\begin{equation}
 	||e_{i+1}|| \leq \bar{\eta} ||e_i||~~\forall e_i \in \bar{\eps}\mc{B}.~~\qed
\end{equation} 
\end{pf*}

Theorem~\ref{thm:GN_convergence} requires that $H_i$ be a sufficiently good approximation of $\nabla_z F(z^*)$ and that the GN subproblems be strongly regular. A sufficient condition for strong regularity is that the QP \eqref{eq:QP} satisfies the LICQ and SSOSC (Theorem~\ref{thm:strong_reg}). In practice, strong regularity can be achieved by a judicious choice of $H_i$. For example, if the Hessian approximation is convex then it is possible to guarantee strong regularity of the subproblems by adding a regularization term, i.e., $H\gets H + \delta I$ for some small $\delta >0$. Then the mapping $H+\delta I + \mc{N}_K$ is strongly monotone, which implies strong regularity \cite[Theorem 2F.6]{dontchev2009implicit}.

\begin{rmk} 
Theorem~\ref{thm:GN_convergence} just requires the Hessian approximation be sufficiently good. One can conceive of useful approximation schemes other than the GN approximation, e.g., $B_i = \nabla_w^2 \phi(w_i)$ when $\phi$ is convex or $B_i = \nabla_w^2 L(\bar{z},x)$ for some fixed $\bar{z}$ sufficiently close to $z^*$.
\end{rmk}

\section{Time-distributed SQP} \label{ss:time_dist_sqp}
In this section we demonstrate that the methods described in Sections~\ref{ss:regularity} and \ref{ss:sequential_quadratic_programming} satisfy the condition of Remark~\ref{rmk:fit_the_framework} and can therefore be used within the framework presented in Section~\ref{ss:optimization}. 

\begin{rmk}
Note that TD-SQP using the GN Hessian approximation and with $\ell = 1$ corresponds to the RTI scheme \cite{diehl2005real}. As such, when specialized to the RTI scheme, Theorems~\ref{thm:iss_proof} and \ref{thm:constraint_satisfaction} are a significant extension of the existing analysis \cite{diehl2005nominal} which does not consider inequality constraints.
\end{rmk}

\textbf{Strong Regularity Assumption:} As detailed in Section~\ref{ss:convex_constraint_sets}, in the presence of convex constraint sets Assumption~\ref{ass:pointwise_strong_reg} can be reduced to the following:
\begin{ass}
The second order sufficient condition \eqref{eq:convex_SOSC} holds at all minimizers in $\Gamma$.
\end{ass} 
As detailed in Section~\ref{ss:NLP}, in the nonlinear inequalities setting, Assumption~\ref{ass:pointwise_strong_reg} can instead be ensured under the following:
\begin{ass}
The linear independence constraint qualification (see Definition~\ref{dfn:LICQ}) and strong second order sufficient condition (see Theorem~\ref{thm:strong_reg}) hold at all minimizers in $\Gamma$.
\end{ass}

\textbf{Algorithm Definition:} Both SQP methods described in Section~\ref{ss:sequential_quadratic_programming} are instances of the following iterative process
\begin{equation}
	H_i(z_{i+1} - z_i) + F(z_i,x) + \mc{N}_K(z_{i+1}) \ni 0,
\end{equation}
for specific choices of $z$, $F$, and $K$. Thus, in both cases the optimization mapping \eqref{eq:opt_def} can be written as
\begin{equation}
	\mc{T}(z,x,i) = (H_i + \mc{N}_K)^{-1} (H_iz - F(z,x)).
\end{equation}

\textbf{Convergence Rate:} If the exact Hessian is used then Theorem~\ref{thm:JN_convergence} applies and the method is at least q-linearly convergent with $q = 2$. If the GN Hessian approximation is used then Theorem~\ref{thm:GN_convergence} applies under some additional assumptions regarding the accuracy of the Hessian approximation, and the method is at least q-linearly convergent with $q = 1$. In both cases the definition of q-linear convergence requires that there be a uniform convergence constant $\eta$ and convergence radius $\veps$ over $\Gamma$. Under the assumption that the functions $\bar{\eta}(x)$ and $\bar{\eps}(x)$ in Theorems~\ref{thm:JN_convergence} and \ref{thm:GN_convergence} are upper and lower semicontinuous, respectively, these can be defined as $\veps = \inf_{x\in\Gamma} \bar{\eps}(x)$ and $\eta = \sup_{x\in\Gamma} \bar{\eta}(x)$. Thus, SQP fits into the framework in Section~\ref{ss:optimization} and can be used for time distributed optimization.

\section{A Numerical Example} \label{ss:numerical_results}

\begin{figure}[htbp]
	\centering
	\includegraphics[width=0.95\columnwidth]{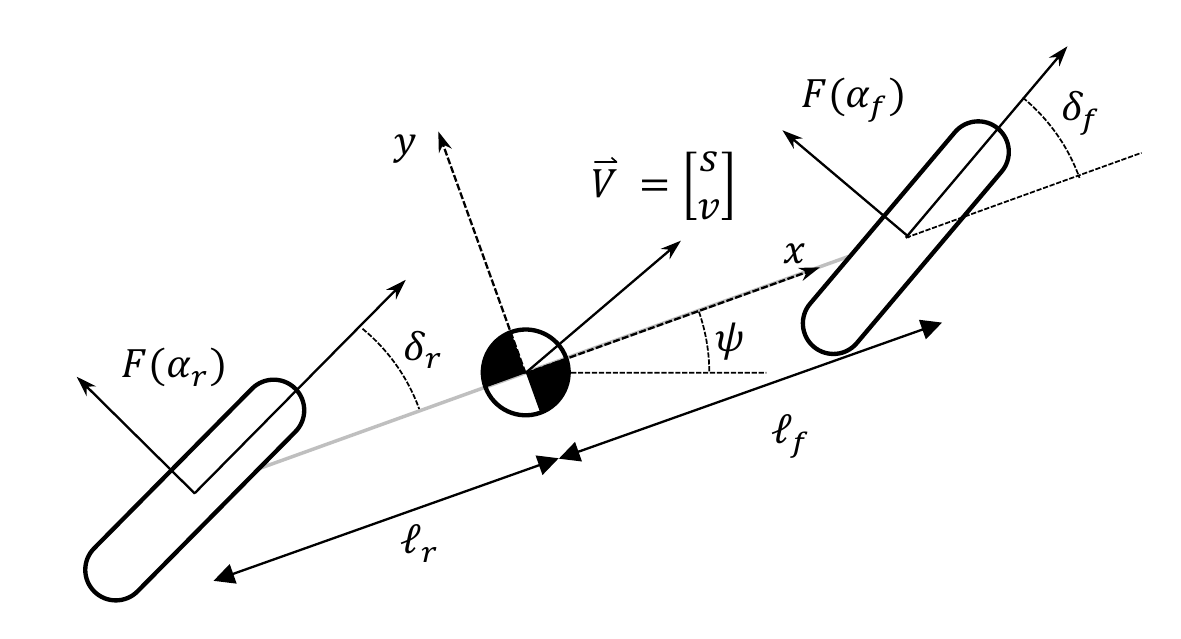}
	\caption{A diagram of the bicycle model}
	\label{fig:bicycle_model}
\end{figure}

Figure~\ref{fig:bicycle_model} illustrates a bicycle model of a sedan. We only consider the lateral portion of the dynamics; the longitudinal velocity $s$ is assumed constant. The states and control inputs are,
\begin{equation}
	x = [y~~\psi~~\nu~~\omega~~\delta_f~~\delta_r], ~~ u = [\dot{\delta}_f~~\dot{\delta}_r],
\end{equation}
where $y$ is the lateral position, $\nu$ is the lateral component of velocity, $\psi$ is the yaw angle, $\omega$ is the yaw rate, $\delta_f$ is the front steering angle, and $\delta_r$ is the rear steering angle.

The equations of motion are
\begin{align*} 
\dot{y} &= s \sin(\psi) + \nu \cos(\psi),\\
\dot{\psi} &= \omega,\\
\dot{\nu} &= -s\omega + \frac{F(\alpha_f) \cos(\delta_f) + F(\alpha_r) \cos(\delta_r) + F_w}{m},\\
\dot{\omega} &= \frac{F(\alpha_f) \cos(\delta_f) \ell_f - F(\alpha_r) \cos(\delta_r) \ell_r}{I_{zz}},\\
\dot{\delta_f} &= \dot{\delta_f},\quad \dot{\delta_r} = \dot{\delta_r},
\end{align*}
where
\begin{align*}
F(\alpha) &= \mu~9.81~m \sin\left(C \arctan(B~\alpha)\right),\\
\alpha_f &= \delta_f - \arctan\left(\frac{\nu + \ell_f \omega}{s}\right),\\
\alpha_r &= \delta_r - \arctan\left(\frac{\nu - \ell_r \omega}{s}\right),\\
&\text{and } F_w = 1/2\rho C_d A |d| d.
\end{align*}
The tire forces are described by a Pacejka model.This model is a modified version of the one presented in \cite{wurts2018collision} and roughly represents a 2017 BMW 740i. The vehicle is disturbed by normally distributed wind gusts $d$ with a mean velocity of $15~m/s$ and standard deviation of $5~m/s$. We obtain a discrete time model using a forward Euler integration scheme with a sampling period of $t_s = 0.04 s$ leading to a discrete time model of the form $x_{k+1} = f_d(x_k,u_k,d_k)$. The model parameters are summarized in Table~\ref{tab:bi_params}\footnote{SI units are used and all angles are in radians unless otherwise noted.}.

\begin{table}[hb]
\caption{Bicycle Model Parameters}
\label{tab:bi_params}
\begin{tabular}{|c|c|c|} \hline
Name & Symbol & Value \\ \hline
Mass & $m$ & $2041~kg$ \\ \hline
Yaw Inertia & $I_{zz}$ & $4964~kgm^2$ \\ \hline
Front, Rear CG distance & $\ell_f,\ell_r$ & $1.56,1.64~m$ \\ \hline
Coefficient of friction & $\mu$ & $0.8$ \\ \hline
Tire parameters & $B,C$ & $12,1.285$ \\ \hline
Lateral Area & $A$ & $7.8~m^2$ \\ \hline
Air Density & $\rho$ & $1.225~kg/m^2$ \\ \hline
Lateral Drag Coefficient & $C_d$ & 1.5 \\ \hline
Longitudinal Velocity & $s$ & $30~m/s$\\ \hline
\end{tabular}
\end{table}

The control objective is to perform a lane change maneuver. This can be achieved by stabilizing the origin which is chosen to coincide with the center of the target lane. The vehicle begins in the neighboring lane at the initial condition $x_0 = \left[-3.7~~0~~0~~0~~0~~0\right]^T$. The OCP is
\begin{subequations}
\begin{align}
\underset{\xi,\mu}{\mathrm{min.}}~~&||\xi_{30}||_{Q_f}^2 + \sum_{i = 0}^{29} ||\xi_i||_Q^2 + ||\mu_i||_R^2,\\
\mathrm{s.t.}~~ &\xi_{i+1} = f_d(\xi_i,\mu_i,0), \quad i= 0,\ldots, 29,\\
& \xi_0 = x(t), ~~ A_f \xi_{30} \leq b_f, \\
& x_{lb} \leq \xi_i \leq x_{ub},  \quad i = 1, \ldots, 30,\\ 
& u_{lb} \leq \mu_i \leq u_{ub}, \quad i=0,\ldots, 29,
\end{align}
\end{subequations}
where $f_d$ is the discrete time model of the sedan. The vehicle is subject to state constraints which keep the vehicle on the road and restrict its yaw and steering angles. The state constraints on $y,\psi,v$ and $\omega$ are softened using $L_1$ exact penalty functions which are implemented using slack variables in order to satisfy our smoothness assumptions. The upper and lower bounds are 
\begin{gather*}
	x_{ub} = [0.4~~7^\circ~~100~~100~~35^\circ~~4^\circ],\\
	x_{lb} = -[4.7~~7^\circ~~100~~100~~35^\circ~~4^\circ],\\
	u_{ub} = [1.2~~0.6], \quad u_{lb} = -[1.2~~0.6],
\end{gather*}
and the weighting matrices are $Q = I_{6\times 6}$, and $R = I_{2\times 2}$. The terminal weight is obtained by solving the discrete time algebraic Riccati equation using the linearization about the origin. The matrices encoding the terminal set, $A_f$ and $b_f$, are computed using the MPT3 toolbox \cite{herceg2013multi}. The natural residual 
\begin{equation} \label{eq:natural_residual}
	\pi(z,x) = ||z - \Pi_K[z - F(z,x)]||,
\end{equation}
is an error bound \cite{pang1997error}, i.e., it upper and lower bounds $||z-z^*(x)||$, where $z^*(x) \in S(x)$, and is commonly used as an easily computable surrogate for the error. We use it throughout this section to measure $||z-z^*(x)||$.

\begin{figure}[htbp]
	\centering
	\includegraphics[width=0.95\columnwidth]{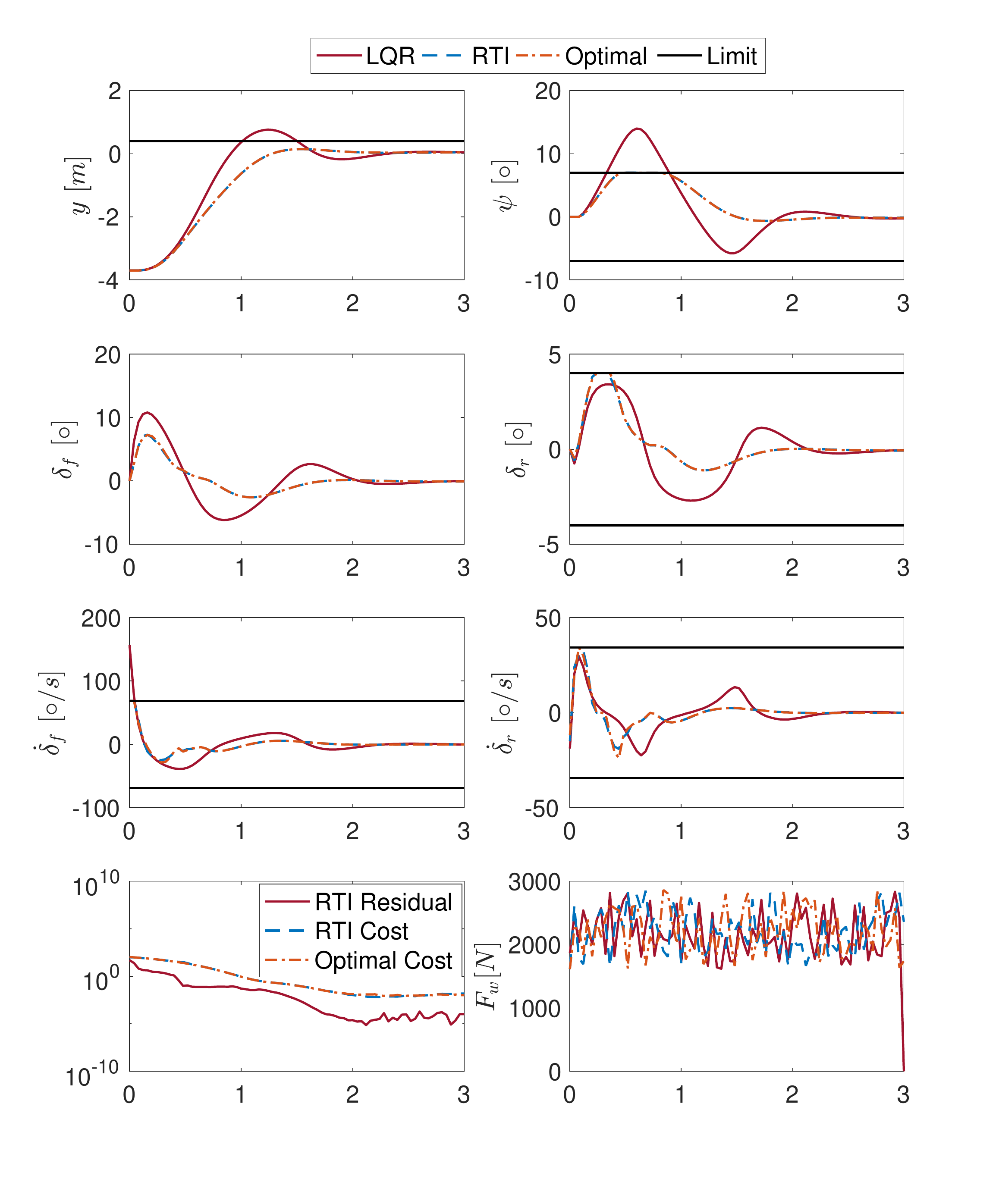}
	\caption{A comparison between an LQR controller, the RTI scheme, and the optimal MPC controller.}
	\label{fig:lqr_rti_opt}
\end{figure}

\begin{figure}[htbp]
	\centering
	\includegraphics[width=0.95\columnwidth]{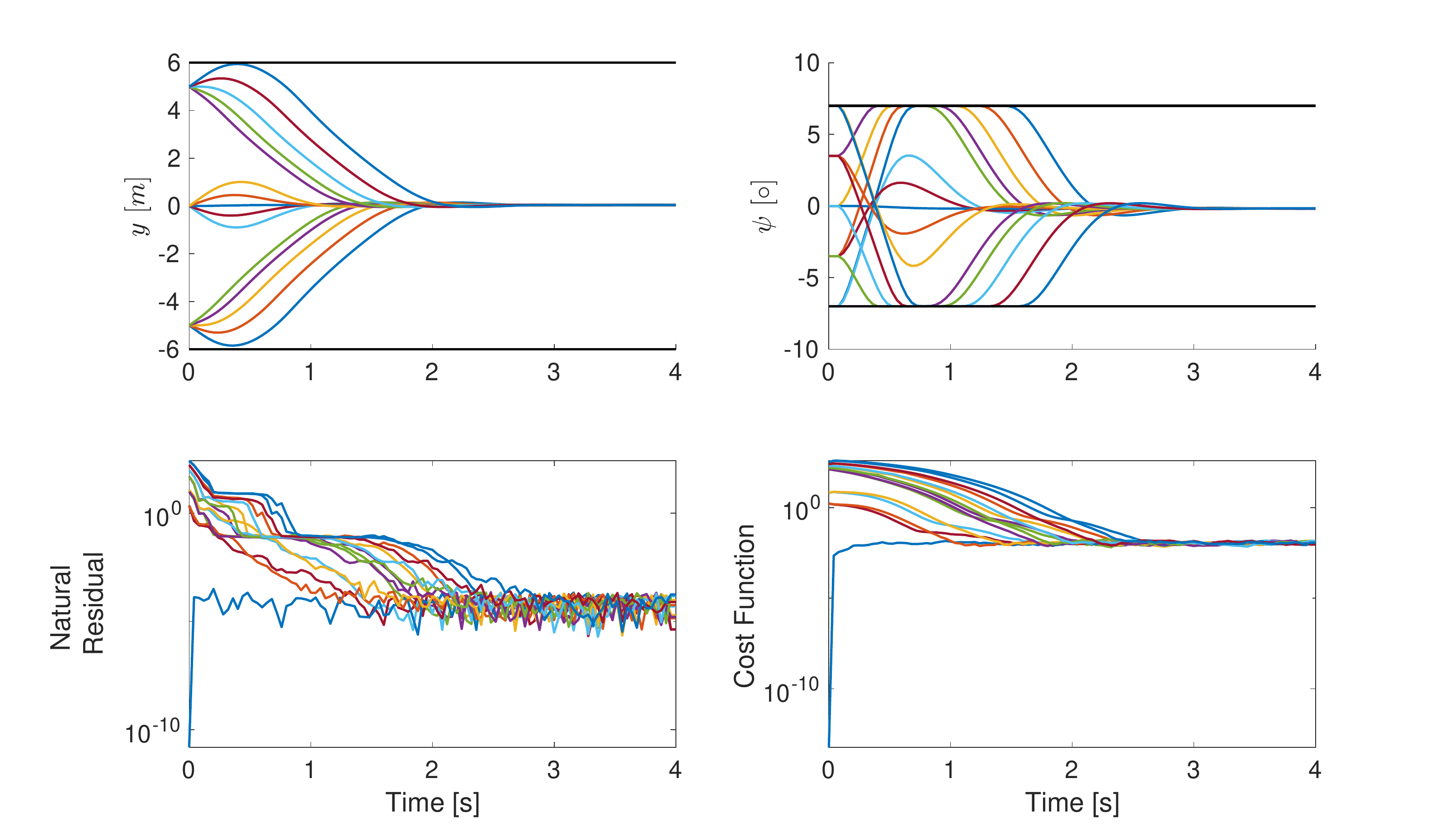}
	\caption{Closed-loop responses for an RTI controller for a variety of initial positions and yaw angles.}
	\label{fig:ROA}
\end{figure}

Figure~\ref{fig:lqr_rti_opt} compares the RTI scheme \cite{diehl2005nominal}, i.e., a TD-SQP scheme using the GN Hessian approximation with $\ell = 1$, with an LQR controller and the optimal MPC feedback law\footnote{All simulations were carried out in MATLAB 2017b on a 2015 Macbook Pro with 16GB of RAM and a 2.8GHz i7 processor. We solved quadratic programs using \texttt{quadprog}. The optimal MPC feedback law was computed using \texttt{fmincon} with default settings. CASADI \cite{andersson2012casadi} was used to compute analytic derivatives which were supplied to the optimization routines.}. The RTI feedback law successfully stabilizes the origin of the plant-optimizer system and outperforms the LQR controller. The state error and the optimization residual both converge to a ball about the origin, demonstrating the expected robustness due to the LISS properties of the combined system (Theorem~\ref{thm:iss_proof}). The closed-loop trajectories generated by the RTI controller are nearly indistinguishable from those from the optimal feedback law but are an order of magnitude cheaper to compute. The RTI scheme took $0.067s$ on average and $0.75s$ in the worst case vs. $0.65s$ and $3.2s$ for the optimal feedback law. Closed-loop responses using the RTI controller for 15 different initial position and yaw angle combinations, with all other states are initialized to zero, are shown in Figure~\ref{fig:ROA}.

\begin{figure}[htbp]
	\centering
	\includegraphics[width=0.95\columnwidth]{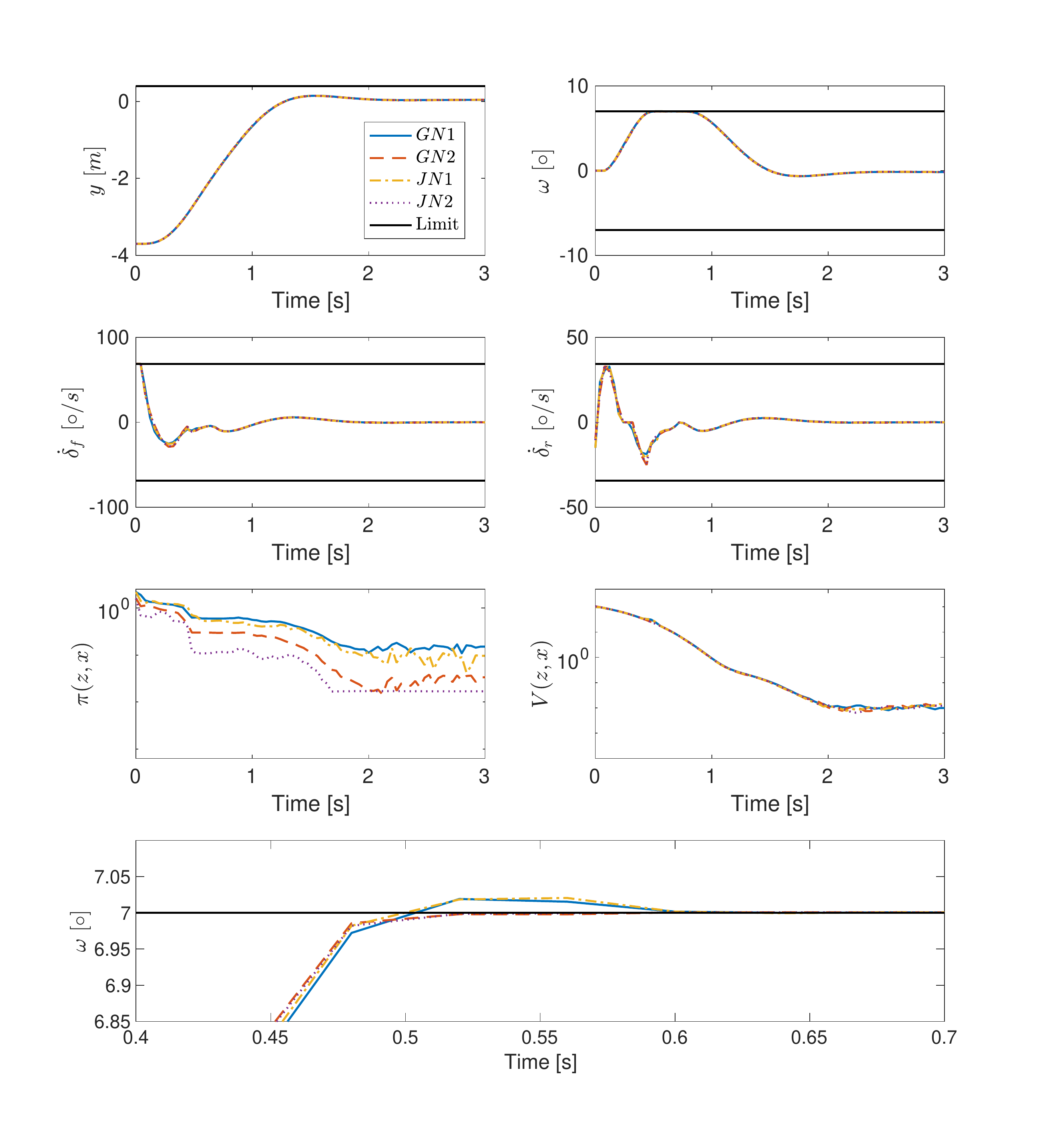}
	\caption{A comparison of TD-SQP controllers implemented using the JN and GN methods.}
	\label{fig:JN_GN}
\end{figure}

Figure~\ref{fig:JN_GN} compares the GN and JN methods with $\ell = 1$ and $\ell = 2$. In the bottom plot of Figure~\ref{fig:JN_GN} note that if $\ell = 2$ iterations are performed, the yaw angle constraint is satisfied exactly, even in the presence of disturbances, as predicted by Theorem~\ref{thm:constraint_satisfaction}. Also, note that the residuals of the computational subsystem converge faster, for a given number of iterations, if the JN method is used instead of the GN method. This is as expected since the convergence rate of the SQP algorithm is faster when the exact Hessian is used.

\section{Conclusions}
In this paper we presented a general framework for the stability analysis of model predictive controllers implemented using time-distributed optimization. When specialized to Sequential Quadratic Programming, our result extends the existing stability analysis of the RTI scheme by explicitly considering inequality constraints, analyzing the effect of performing additional SQP iterations, considering a wider class of Hessian approximations, and proving local input-to-state stability of the closed-loop system. Future work includes analyzing the effect of the sampling rate, applying our framework to globalized SQP methods, and developing numerical methods for estimating the the asymptotic gain functions used in the analysis.

\bibliography{iss_mpc}

\end{document}